\documentclass{article}

\usepackage{a4}
\usepackage{verbatim}
\usepackage{amsmath}
\usepackage{amssymb}

\newenvironment{prf}{\par\noindent{\bf Proof.}}{\par\rightline{$\Box$}}
\def\noprf{\par\rightline{$\Box$}}

\numberwithin{equation}{section}
\newtheorem{theo}[equation]{Theorem}
\newtheorem{lem}[equation]{Lemma}
\newtheorem{cor}[equation]{Corollary}
\newtheorem{con}[equation]{Conjecture}

\newtheorem{defin}{Definition}

\newtheorem{ex}[equation]{Example}
\newtheorem{rem}[equation]{Remark}
\newtheorem{prop}[equation]{Proposition}

\def\ideal{\vartriangleleft}

\def\uD{\mathrm{D}}
\def\ud{\mathrm{d}}
\def\sspan{\operatorname{span}}
\def\pd{\partial}
\def\Sp{\mathbf{Sp}_{\F}(V)}

\def\Gl{\mathbf{Gl}}
\def\exp{\operatorname{exp}}
\def\Ad{\operatorname{Ad}}
\def\Id{\operatorname{Id}}
\def\I{\mathcal{I}}
\def\J{\mathcal{J}}
\def\ccS{\mathcal{S}}

\def\Sym{\operatorname{Sym}}

\def\fii{\varphi} 
\def\lam{\lambda}  
\def\m{\mu}     

\def\mono{\hookrightarrow}

\def\ad{\mathfrak{ad}}
\def\gotg{\mathfrak{g}}
\def\goth{\mathfrak{h}}
\def\gotw{\mathfrak{w}}
\def\gotp{\mathfrak{p}}
\def\gots{\mathfrak{s}}
\def\gotsp{\mathfrak{sp}}
\def\gotsl{\mathfrak{sl}}
\def\gotso{\mathfrak{so}}
\def\gotgl{\mathfrak{gl}}

\def\gotn{\mathfrak{n}}
\def\gota{\mathfrak{a}}
\def\gotb{\mathfrak{b}}
\def\gote{\mathfrak{e}}

\def\half{\frac{1}{2}}

\def\C{{\mathbb C}}
\def\F{{\mathbb C}}

\def\N{{\mathbb N}}

\def\R{{\mathbb R}}
\def\Z{{\mathbb Z}}

\def\P{{\mathbb P}}

\def\ra{\rightarrow}
\def\lra{\longrightarrow}

\def\osum{\oplus}
\def\tensor{\otimes}

\def\izo{\simeq}

\def\ccF{{\cal F}}

\def\ccL{{\cal L}}

\def\ccO{{\cal O}}

\def\ccS{{\cal S}}

\def\diag{\operatorname{diag}}

\newcommand{\mr}[1]{\rightarrow}

\usepackage[all]{xy}
\usepackage[dvips]{graphicx}

\author{Jaros\l{}aw Buczy\'nski}

\title{Properties of legendrian subvarieties of projective space}

\begin{document}

\maketitle

\begin{abstract}
I prove that every smooth legendrian variety generated by quadrics is a homogeneous variety and further I give 
a list of all such legendrian varieties. A review of the subject is included, illustrated by examples.
 Another result is that no complete 
intersection is a legendrian variety. 
\end{abstract}

\tableofcontents


\section{Introduction}

The main result of this article is the full classification of smooth legendrian varieties with ideal generated
by quadratic polynomials. The key fact for this is that such varieties are homogeneous. 
This is done by constructing a subgroup of $\Sp$ out of the quadratic part of the ideal.
Next the subgroup acts on the legendrian variety and if the ideal is generated by it's quadratic part,
then the action is transitive on smooth points.

Another result says that if a legendrian variety is a complete intersection then it is either a linear space
 or it's singular locus is of codimension 1.

The article is meant to be elementary and contains a lot of well known statements. An expert reader should skip to
the subsection \ref{section_proj_geometry} for the second result mentioned or to the sections \ref{algebry_lie},
\ref{rozdzial_rozm_gen_kwadryki} and \ref{quadrics} for the first result. 

The structure of the paper is as follows:

In the section \ref{rozdzial_poczatkowy} I define the main notions of the paper.

In the section \ref{curves} I explain what is known about the legendrian curves. This is not at all essential to the 
rest of paper (except possibly the example \ref{przyklad_skr_kubika}), but it might give an idea of the subject. 

The section \ref{general_remarks} analyzes some well known notions (such as conormal bundle, Atiyah extension...) 
in the case of legendrian varieties and finds relations between them finally leading to the proof of the theorem about 
legendrian complete intersections.

Next, the section \ref{algebry_lie} recalls the Poisson bracket on a polynomial ring. 
Further some subalgebras (of the Poisson structure)
 which at the same time are ideals (with respect to the standard multiplication)
are related to legendrian varieties. 
The restriction to the quadratic part of the ring calls for finite dimensional Lie algebra $\gotsp_{2n}$
 and for group $\Sp$. 
This is already the section \ref{rozdzial_rozm_gen_kwadryki}, 
which takes reader to the proof of the key theorem of the article 
- the theorem \ref{rozmaitosci_jednorodne}. 

The section \ref{rozdzial_przyklad} explains motivation for the study of the legendrian varieties.
Further, I describe a fabulous family of examples - the subadjoint varieties. 
They turn up to have yet another description - they are the legendrian varieties generated by quadrics.
The explicit formulations about the examples where computed with a help of Magma program.

The section \ref{representation_theory} revises some representation theory.

Finally, in the section \ref{quadrics} I use the theorem  \ref{rozmaitosci_jednorodne} and the techniques of
the Lie theory to prove the main result of the article. Unfortunately it is quite technical.

\subsection{Acknowledgements}

The heart of the article is just a revised, corrected and extended translation of my Master Degree dissertation \cite{mgr},
supervised by Jaros\l{}aw Wi\'sniewski. 
He predicted most of the results and was a tremendous  support during the whole, long process of writing the article.

I am grateful to a number of people for their help in understanding the Lie theory: 
Colin Ingalls, Mariusz Koras, Adrian Langer, Miko\l{}aj Rotkiewicz, Dmitriy Rumynin. 

Partially this article was created during my stay at Warwick University, as a Marie Curie fellow.
I would like to thank Miles Reid for his invitation and hospitality and Mark Gross for all his help and support.
Special thanks to Gavin Brown for teaching me how to use Magma. 

Finally, I would like to mention a lot of enlightening discussions with Luis \mbox{Eduardo} Sola Conde - THANKS!

And the very last acknowledgement to my wife Weronika for her patience, empathy and tolerance.

\section{Basic definitions}
\label{rozdzial_poczatkowy}

In all the definitions below and through all the article I assume that $V$ is a finite dimensional vector space over 
a field $\C$  of complex numbers. The only exceptions are example \ref{example_real} and remark \ref{real_remark}, where 
the ground field is $\R$. 

\begin{defin}
 A \textbf{symplectic form} on the vector space $V$ is a map $\omega : V \times V \ra \F$ which is:\\
a) bilinear,\\
b) skew-symmetric and\\
c) non-degenerated (meaning $\ \forall_{v\in V} \ \exists_{w\in V} \ \omega(v,w) \ne 0$).
\end{defin}

Distinguishing a basis of $V$ gives a correspondence between bilinear forms and matrices.
 From now on I will identify symplectic form with its matrix.

If there exists a symplectic form on the vector space $V$ then its dimension is even. Hence from now on I will assume that $\dim V= 2n$ for an integer $n$.

\begin{defin}
Let $\omega$ be a symplectic form on $V$. A linear subspace $W \subset V$ is a \textbf{Lagrangian subspace}, if it is a maximal subspace such that $\omega |_W
\equiv 0$ (so in particular $\dim W = \half \dim V = n$). I will denote a Lagrangian subspace via $W\subset_L V$.
\end{defin}

Given an algebraic variety $X$ I denote by $X_0$ its smooth locus. If $X \subset \P(V)$ is a projective variety, 
I denote by $\hat{X} \subset V$ the affine cone over $X$.
For convenience by a slight abuse of this notation let $\hat{X}_0 \subset V$ denote smooth points of the cone over 
$X$ except possibly $0 \in V$, i.e. 
$\hat{X}_0:= (\hat{X})_0 \backslash \{0\}$.

\begin{defin}
A subvariety $X \subset \P(V)$ is \textbf{legendrian} (denoted by $X \subset_l \P(V)$) if for each smooth point of its affine cone the tangent space at this point is 
Lagrangian. In other words:
$$X\subset_l \P(V) \iff \forall_{x \in \hat{X}_0} T_x \hat{X}_0 \subset_L T_x V = V$$
\end{defin}

\begin{rem}
The dimension of a legendrian subvariety in $\P(V)= \P^{2n-1}$ is equal to $n-1$.
\end{rem}

\begin{ex}
If $n=1$ then any point of $\P^1$ is a legendrian subvariety.
\end{ex}

\begin{rem} \label{nieprzywiedlne}
A variety is legendrian if and only if each of its irreducible components is legendrian.
\end{rem}
\nopagebreak
\noprf
\begin{cor}
$X \subset \P^1$ is a legendrian subvariety if and only if $X= \{p_1, \ldots p_k\}$ for some $p_1,\ldots, p_k \in \P^1$.
\end{cor}
\noprf
\begin{rem} \label{uwaga_o_otw_podzb}
For an irreducible variety $X$ to be legendrian it is enough that the condition from the definition is satisfied 
for every point in an open subset $U\subset \hat{X}_0$.
More formally: 
$$X\subset_l \P(V) \ \iff \ \exists_{U\subset V} \ \forall_{x \in U} \ T_x \hat{X}_0 \subset_L V$$
\end{rem}
\noprf

\begin{ex}\label{przyklad_liniowy}
Let $W \subset V$ be a linear subspace then 
$$ \P(W) \subset_l \P(V) \iff  W \subset_L V $$
\end{ex}

\begin{theo}\label{o_rozm_zdegenerowanych}
Let $X \subset_l \P(V)$. Then the following conditions are equivalent:
\begin{itemize}
\item[(i)]
$X$ is degenerated (i.e. is contained in a hyperplane).
\item[(ii)]
There exists a linear subspace $W' \subset V$ of codimension 2 such that $\omega|_{W'}$ is a symplectic form on $W'$, $X' = \P(W') \cap X$ is a legendrian subvariety 
in $\P(W')$ and $X$ is a cone over $X'$.
\item[(iii)]
$X$ is a cone over some variety $X'$.
\end{itemize}
\end{theo}

\begin{prf}
First prove that (i) implies (ii). Let $X \subset \P(W)$ for linear subspace $W \subset V$ of codimension 1. Then W contains one dimensional linear subspace $Q$
such that for all $q\in Q$ and $w\in W$ one has $\omega(q,w)=0$. Hence for each point $x \in \hat{X}_0$ line $Q$ is contained in $T_x \hat{X}_0$. So $Q\in X$ and $X$ consists of 
lines through Q. Now let $W'$ be any subspace of $W$ not containing $Q$. Verifying that $W'$ satisfies (ii) I leave to reader.\\

Since the implication (ii) $\Longrightarrow$ (iii) is obvious, it remains to prove the implication from (iii) to (i). Let $Q$ be the vertex of the cone. Then $Q$ is a line in $V$, which is contained
in $T_x \hat {X}_0$ for all $x \in \hat{X}_0$. Let $W$ be the linear subspace of $V$ perpendicular (with respect to $\omega$) to $Q$ 
(so $W=\{w \in V | \forall_{q \in Q} \  \omega (w,q) = 0\}$). Since $\omega |_{T_x \hat{X}_0} \equiv 0$ it follows that $T_x \hat{X}_0 \subset W$ for all $x \in \hat{X}_0$,
so $X \subset \P(W)$
\end{prf}

\section{Legendrian curves}\label{curves}

\begin{ex}\label{przyklad_skr_kubika}
(twisted cubic) Consider the Veronese embedding of degree 3: \mbox{$\P^1 \stackrel{\fii}{\mono} \P^3$}, defined by $\fii(\lam : \m)= (\lam^3 : \lam^2 \m : \lam \m^2 : \m^3)$.
If the symplectic form $\omega$ is defined by the matrix
\begin{displaymath}
\left(\begin{array}{rrrr} 0&0&0&-1\\
                    0&0&3&0\\
                    0&-3&0&0\\
                    1&0&0&0 \end{array} \right)
\end{displaymath}
then the curve $X=\fii(\P^1)$ is legendrian.
\end{ex}
\begin{prf}
Fix a point $x \in \hat{X}_0$, so that $x=(\lam_0^3 , \lam_0^2 \m_0,\lam_0 \m_0^2 , \m_0^3)$ for some \mbox{$(\lam_0,\m_0) \in \C^2\backslash \{0,0\}$}.
Tangent space to $\hat{X}_0$ at $x$ is spanned by vectors
$$v_1:=\frac{\pd\fii}{\pd\lam}(\lam_0,\m_0) = (3\lam_0^2, 2\lam_0 \m_0, \m_0^2,0)$$
\nopagebreak
$$v_2:=\frac{\pd\fii}{\pd\m}(\lam_0,\m_0) = (0, \lam_0^2,2 \lam_0 \m_0, 3 \m_0^2)$$
 Now calculate:
\begin{displaymath}
\omega(v_1,v_2) \; = \; -3 \m_0^2 \cdot \lam_0^2 \; + \; 3 \cdot
(2\lam_0 \m_0 ) \cdot (2\lam_0 \m_0) \; - \; 3\cdot
(3\lam_0^2)\cdot(\m_0^2) \; = \; 0.
\end{displaymath}
\textit{Remark: It is an easy exercise to show that if any other symplectic form $\tilde{\omega}$ satisfies $\tilde{\omega}(v_1,v_2)=0$ for all $\lam_0,\m_0 \in \C^2
\backslash \{0\}$, then it equal to $\omega$ (up to proportionality), so $\tilde{\omega}=\alpha \omega$ for some $\alpha \in \C^*$.}

So by bilinearity and skew-symmetry of $\omega$ we obtain that $\omega|_{T_x \hat{X}_0} \equiv 0$, so indeed $X \subset \P^3$ is a legendrian subvariety.
\end{prf}

\begin{theo}\label{classification_n=2}
Let $X \subset \P^3$ be an irreducible legendrian curve. If X is generated by the polynomials of degree at most 2,
then $X\izo \P^1$ and $X$ is either a line (see the example \ref{przyklad_liniowy}) or a twisted cubic (see the example \ref{przyklad_skr_kubika}).
\end{theo}

\begin{prf} I examine case by case. First notice, that for degenerated (i.e contained in a plane) curves the situation 
is clear: By the theorem 
\ref{o_rozm_zdegenerowanych} there exist $X'$ such that $X$ is a cone over $X'$. Since $X$ is irreducible, it follows 
that $X'$ is a point and so $X$ is a line.

Second, let $d:=\deg X$ and now notice, that the degree of a curve in $\P^3$ generated by quadrics is bounded: it is less or equal to 4. 
Moreover, if $d=1$ or $d=2$, then the curve $X$ is surely degenerated (see for example \cite[IV.3.11.1]{hartshorn}).
 So assume $X$ is not degenerated and so either 
$d=3$ or $d=4$. For $d=3$ the only non-degenerated curve in $\P^3$ is the twisted cubic, so it remains to exclude 
the case $d=4$.

But $d=4$ is the case of complete intersection. If $X$ is a singular complete intersection of two quadrics, then it is 
a rational curve which for some choice of coordinates is parametrised by:
$$\P^1 \ni (t:s) \mapsto (t^4: t^3s+ts^3: t^2s^2: s^4) \in \P^3$$
(geometrically, it is an intersection of a quadratic cone and a sphere, such that the sphere passes through the vertex of the 
cone)\footnote{I would like to thank Jorge Caravantes for this enlightening remark. Also I am really grateful to 
Luis Eduardo Sola Conde for his help in calculating explicitly the parametrisation.}. Now it is really 
easy to compute few tangent spaces to the cone over the curve and verify that they cannot be altogether Lagrangian
 with respect to any symplectic form.

The last case is a smooth complete intersection. Again a coordinate system might be chosen, so that the quadrics are:
$$q_1:=x_0^2 + x_1^2 + x_2^2 \textrm{ and } q_2:= x_1^2 + \lam x_1 x_2 + x_2^2 +x_3 ^2$$
for some $\lam \in \C^*$ (a bit different choice of coordinates might be found in \cite[prop. 2.1(d)]{miles}). Again it is easy 
to take few tangent spaces to the cone over such a curve and verify that they cannot be altogether Lagrangian with respect to any
symplectic form. 
\end{prf}

\subsection{Rational curves}

\begin{ex}\label{przyklad_krzywe_legendrowskie}
Let $k$ and $l$ be coprime natural numbers. Consider a map $\fii_{k,l} :\P^1 \lra \P^3$ defined by
$$\fii_{k,l}(\lam:\m)= (\lam^{k+l}: \frac{k-l}{k+l} \m^{k+l}: \lam^l \m^k:\lam^k \m^l).$$
By $X_{k,l}$ denote the image $\fii_{k,l}(\P^1)$. If the symplectic form is given by the matrix 
$$\begin{pmatrix}
  0 & 1 & 0 & 0 \\
  -1& 0 & 0 & 0 \\
  0 & 0 & 0 & 1 \\
  0 & 0 & -1& 0
\end{pmatrix}$$
then $X_{k,l}$ is a legendrian curve. In particular for $k=2$, $l=1$ we get a twisted cubic and for $k=1$, $l=0$ -- line.
\end{ex}

This might be verified in exactly the same way as for the twisted cubic. But instead of repeating the argument I will describe rough classification of rational legendrian 
curves.\\

Consider a rational curve $X\subset \P^3$, which is the closure of the image of map 
$$
\begin{array}{rcl}
\C &\lra&\P^3\\
t &\longmapsto& \big( 1:f_1(t):f_2(t):f_3(t) \big)
\end{array}
$$
for some polynomials (not necessarily homogeneous) $f_1,f_2,f_3$. The cone $\hat{X}$ is the closure of the set of points of the form $\big(\alpha, \alpha f_1(t),
\alpha f_2(t), \alpha f_3(t)\big)$ for some $\alpha \in \C$, $t\in \C$. The tangent space at such a point (if only it is smooth) is spanned by two vectors:\\
\indent $v_1 = \big(1, f_1(t), f_2(t), f_3(t)\big)$ (direction tangent to generator of the cone) and\\
\indent $v_2 = \big(0, \frac{\ud f_1}{\ud t} (t), \frac{\ud f_2}{\ud t} (t), \frac{\ud f_3}{\ud t} (t)\big)$ (direction tangent to the base of the cone).\\
By the remark \ref{uwaga_o_otw_podzb}, $X$ is legendrian if and only if for all $t\in \C$ I have $\omega(v_1,v_2)=0$. Assume that $\omega$ is as in the example 
\ref{przyklad_krzywe_legendrowskie} and compute:
$$\omega(v_1,v_2) = \frac{\ud f_1}{\ud t} (t) + f_2(t)\frac{\ud f_3}{\ud t} (t) - f_3(t)\frac{\ud f_2}{\ud t} (t)$$

Therefore $X$ is a legendrian curve if and only if polynomials $f_1, f_2, f_3$ satisfy a differential equation:

\begin{equation}\label{rownanie_rozn_dla_krzyw_wym}
 \dot{f_1} = \dot{f_2} f_3 - \dot{f_3} f_2
\end{equation}

It is clear that there is a lot of solutions of this equation: given any polynomials $f_2$ and $f_3$ we can find $f_1$ satisfying \ref{rownanie_rozn_dla_krzyw_wym}. In 
particular for $f_2:=t^k$, $f_3:= t^l$, the resulting curve is $X_{k,l}$ from the example \ref{przyklad_krzywe_legendrowskie}.\\

Exactly the same argument can be applied to rational functions instead of polynomials. And therefore I get:
\begin{theo}\label{wymierne_krzywe_legendrowskie}
A rational curve $X\subset \P^3$ parametrised by 
$$\P^1\ni \big(\lam:\m \big) \lra \big(\fii_0 (\lam,\m):\fii_1 (\lam,\m):\fii_2 (\lam,\m):\fii_3 (\lam,\m)\big) \in \P^3$$
is legendrian if and only if either $\fii_0 \equiv 0$ and $X$ is a legendrian line or the rational functions 
$f_i(t):= \frac{\fii_i(t,1)}{\fii_0(t,1)}$ satisfy the equation \ref{rownanie_rozn_dla_krzyw_wym}.
\end{theo}
\noprf

\subsection{Other curves}

A very similar statement as the theorem \ref{wymierne_krzywe_legendrowskie} is true for all the other curves - 
see \cite[thm F]{bryant}, This finally leads to proof that every smooth irreducible curve admits a legendrian embedding in
$\P^3$ \cite[thm G]{bryant}. However, note that this is not true for higher dimensions, although every variety is birational to a 
(usually singular) variety. For more details consult \cite{landsbergmanivel04}.


\section{General remarks}\label{general_remarks}
In this section I will explain few more conditions equivalent to the definition of a legendrian subvariety and also I shall present some new examples.

\subsection{Tangent and conormal bundles}
The definition of legendrian variety given in the section \ref{rozdzial_poczatkowy}
 is using the tangent space of $\hat{X}$. Sometimes it is more convenient to have 
a dual condition, involving conormal bundle, since it is an image of the ideal of the variety.
 So here in this subsection I present an equivalent statement for $X$ to be 
legendrian.

Let $\omega$ be a symplectic form on a vector space $V$. I start pointing out that the form $\omega$ can be expressed as 
 an isomorphism of vector spaces: $V$ and its dual $V^*$:\\
\begin{equation} \label{izomorfizm_phi}
\begin{array}{cccl}
 \phi :&V&\stackrel{\izo}{\lra}&V^*\\
       &v&\longmapsto          &\omega(v,\ \cdot \ )
\end{array}
\end{equation}

Let me denote by $\omega' : V^* \times V^* \ra \C$ the pullback of $\omega$ using $\phi^{-1}$:
\begin{displaymath}
\omega':= (\phi^{-1})^* \omega = \omega \big( \phi^{-1}(\cdot),
\phi^{-1}(\cdot) \big)
\end{displaymath}

\begin{theo} \label{tangent_and_conormal}
A subvariety $X \subset \P(V)$ is legendrian if and only if for each smooth point of its affine cone the conormal space 
at this point is Lagrangian with respect to form $\omega'$.
\end{theo}
\begin{prf}
Consider the following diagram:

\begin{equation} \label{diagram_wiazek}
\xymatrix{ 0 \ar@{->}[r] & T\hat{X}_0 \ar@{->}[r]^i
\ar@{.>}[d]^{\psi} &TV|_{\hat{X}_0} \ar@{->}[r]^p
\ar@{->}[d]^{\phi}_{\izo}& N_{\hat{X}_0 \slash V} \ar@{->}[r]&0\\
 0 \ar@{->}[r] & N^*_{\hat{X}_0 \slash V} \ar@{->}[r]^j
&T^*V|_{\hat{X}_0} \ar@{->}[r]^q& T^* \hat{X}_0 \ar@{->}[r]&0 }
\end{equation}

The rows of the above diagram are exact. $TV|_{\hat{X}_0}$ is a trivial vector bundle over $\hat{X}_0$ with fibre V and similarly $T^*V|_{\hat{X}_0}$ is a trivial
vector bundle with fibre $V^*$. So $\phi$ is the isomorphism in each fibre defined as in the equation \eqref{izomorfizm_phi}.
Below I construct such $\psi$ that the diagram is commutative.\\

Fix a point $x\in \hat{X}_0$. First assume $X$ is legendrian. I claim that the composition \mbox{$q \circ \phi \circ i$} is zero. Indeed, let \mbox{$v \in T_x \hat{X}_0$}
be any vector and $\alpha \in T^*_x \hat{X}_0$ be its image under that composition:
\begin{displaymath}
\alpha:= q \circ \phi \circ i (v),
\end{displaymath}
then for any $w \in T_x \hat{X}_0$ by lagrangianity of $X$ I have $\alpha(w)= \omega(v,w)=0$.\\

So $q \circ (\phi \circ i) = 0 $ and by the definition of kernel of $q$ there exists  exactly one 
homomorphism $\psi: T\hat{X}_0 \lra N^*_{\hat{X}_0 \slash V}$, such that 
$j\circ \psi = \phi \circ i$. Since $\phi \circ i$ is a monomorphism it follows that $\psi$ is a monomorphism as well.
 But both rank of $T\hat{X}_0$ and rank of 
$N^*_{\hat{X}_0 \slash V}$ are equal to $n$ and so $\psi$ is an isomorphism. Hence $j= \phi \circ i \circ (\psi^{-1})$ and 
\begin{equation}
\omega'|_{N^*_{\hat{X}_0 \slash V}} = j^* \omega' = (\psi^{-1})^*
i^* \phi^* \omega' = (\psi^{-1})^* i^* \omega = (\psi^{-1})^*
\omega|_{T\hat{X}_0} \equiv 0.
\end{equation}

Exactly the same argument applied to the diagram
\begin{displaymath}
\xymatrix{ 0 \ar@{->}[r] & N^*_{\hat{X}_0 \slash V} \ar@{->}[r]^j
\ar@{.>}[d] &T^*V|_{\hat{X}_0} \ar@{->}[r]^q
\ar@{->}[d]^{\phi^{-1}}_{\izo}& T^* \hat{X}_0 \ar@{->}[r]&0 \\
0 \ar@{->}[r] & T\hat{X}_0 \ar@{->}[r]^i &TV|_{\hat{X}_0}
\ar@{->}[r]^p & N_{\hat{X}_0 \slash V} \ar@{->}[r]&0}
\end{displaymath}
shows that if $\dim \hat{X} = 2n - n = n$ and $\omega'|_{N^*_{\hat{X}_0 \slash V}} \equiv 0$, then 
$\omega|_{T_{\hat{X}_0}} \equiv 0$.
\end{prf}

Actually the above proof shows a little more:
\begin{cor}\label{completing_the_diagram}
$X \subset \P(V)$ is legendrian if and only if there exist an isomorphism $\psi$ completing the diagram \eqref{diagram_wiazek}.
\end{cor}
\noprf

\subsection{Projective geometry and complete intersections} \label{section_proj_geometry}

This subsection is not essential to the rest of the paper, but rather explains legendrian subvarieties in a wider context. First I will present an 
equivalent condition for a variety to be legendrian in terms of projective geometry only. Next I prove the theorem \ref{o_zupelnym_przecieciu}, which in particular
 states, that no smooth complete intersection can be a legendrian subvariety.\\ 

Let $X\subset \P(V)$ be any projective variety. Notice that $\hat{X}_0$ is a $\C^*$ principal bundle over $X_0$ which is just a line bundle $\ccO_{X_0}(-1)$ 
with the zero section removed. So there is a natural action of $\C^*$ on $\hat{X}_0$ and on all the vector bundles in the diagram \eqref{diagram_wiazek}. 
In particular, the action has weight $+1$ on $V$ and weight $-1$ on its dual $V^*$,
 i.e. for $t \in \C^*$:
$$t \cdot v = tv \textrm{ for } v \in V \textrm{ and}$$
$$t \cdot \alpha = t^{-1}\alpha \textrm{ for } \alpha \in V^*.$$
Let $\pi: \hat{X}_0 \lra X_0$ denote the bundle projection and instead of vector bundles consider locally free sheaves, so that it makes sense to push them forward by $\pi$ 
and hence I get a diagram:
$$\xymatrix{ 0 \ar@{->}[r] &
\pi_*(T\hat{X}_0) \ar@{->}[r]
 &\pi_*(V\tensor \ccO_{\hat{X}_0}) \ar@{->}[r]
\ar@{->}[d]^{\pi_*(\phi)}_{\izo}& \pi_*(N_{\hat{X}_0 \slash V}) \ar@{->}[r]&0\\
 0 \ar@{->}[r] & \pi_*(N^*_{\hat{X}_0 \slash V}) \ar@{->}[r]
&\pi_*(V^* \tensor \ccO_{\hat{X}_0}) \ar@{->}[r]& \pi_*(T^*
\hat{X}_0) \ar@{->}[r]&0 }$$

All these sheaves again admit the action of $\C^*$ and so each of them decomposes to a direct sum of subsheaves with 
a fixed gradation of this action. More precisely:
$$\ccF = \bigoplus_{i=-\infty}^{\infty} \ccF_i$$
where $\ccF$ is any of the sheaves of the above diagram and $\ccF_i$ is such a subsheaf of $\ccF$ that $\C^*$ acts on
 it with weight $i$, so for an open subset 
$U \subset X_0$ and $t \in \C^*$:
$$t\cdot s = t^i s \  \textrm{ for } \ s \in \ccF_i(U).$$
The isomorphism 
$$\pi_*(\phi): \pi_*(V \tensor \ccO_{\hat{X}_0}) \lra \pi_*(V^* \tensor
\ccO_{\hat{X}_0})$$
switches the gradation by $-2$. So if I restrict in the upper row to the gradation $+1$ and in the lower row to 
the gradation $-1$, I get:
\begin{equation}\label{diagaram_popchniety}
\xymatrix{
 0  \ar@{->}[r]  &  \ccL^*(-1) \ar@{->}[r]                 &  V \tensor \ccO_{X_0} \ar@{->}[r] \ar@{->}[d]_{\izo}  &  N_{X_0 \slash \P(V)}(-1) \ar@{->}[r]  &  0\\
 0  \ar@{->}[r]  &  N^*_{X_0 \slash \P(V)}(1) \ar@{->}[r]  &  V^* \tensor \ccO_{X_0} \ar@{->}[r]                   &  \ccL(1) \ar@{->}[r]                   &  0 }
\end{equation}
where $\ccL$ is the Atiyah extension 
$$ 0 \ra T^*X_0 \ra \ccL \ra \ccO_{X_0} \ra 0$$
corresponding to the Chern class $c_1(\ccO(1))$. The above argument is a mimic of \cite[sect. 2.1]{wisniewski}, but 
the original reference is \cite{atiyah}).

\begin{cor}
Rows of the diagram \eqref{diagaram_popchniety} are isomorphic if and only if rows of the diagram \eqref{diagram_wiazek} are isomorphic. So by the corollary 
\ref{completing_the_diagram} they are isomorphic if and only if the variety $X$ is legendrian.
\end{cor}
\noprf

So I have explained an equivalent definition of legendrian variety in terms of projective geometry only. 
This is actually very close to yet another definition involving the contact form on $\P^{2n-1}$ - the one that I should 
have started from, since it is the original definition, which can be generalised to subvarieties of contact manifolds. 
But in the context of this paper it is more convenient to work with the affine definition given in section 
\ref{rozdzial_poczatkowy} and the contact form on $\P^{2n-1}$ is of little use here. A reader interested in the topic should 
have a look at \cite{landsbergmanivel04}, \cite{wisniewski}, \cite{kebekus1} and \cite{kebekus2} and many other works on the
topic.\\

The construction of the bundle $\ccL$ (and its naturality) 
assures that it is a non-trivial extension (i.e. non-splitting) if only it is restricted to a projective subvariety:
\begin{lem}\label{o_wiazce_L}
If $X\subset_l\P(V)$ and $Y\subset X_0$ is projective (i.e. $Y$ closed in $\P(V)$) with $\dim Y > 0$, then the vector
 bundle $\ccL|_Y$ on $Y$ is a non-trivial extension:
$$ 0 \ra T^*X|_Y \ra \ccL|_Y \ra \ccO_Y \ra 0.$$
Moreover since $\ccL \izo N_{X_0 \slash \P(V)} \tensor \ccO_{X_0}(-2)$, obviously 
$$\ccL|_Y \izo (N_{X_0 \slash \P(V)})|_Y \tensor \ccO_{Y}(-2)$$
\end{lem}
\noprf

\begin{theo} \label{o_zupelnym_przecieciu}
Assume that $X\subset_l \P(V)$ is an irreducible normal complete intersection. Then $X$ is a linear subspace.
\end{theo}
\begin{prf}
Assume that $X$ is a complete intersection of hypersurfaces of degrees $d_1, d_2, \ldots, d_n$. Firstly, I am going to show,
 that at 
least one of the degrees $d_i$ is equal to 1. It is well known, that for a complete intersection:
$$N_{X_0 \slash \P(V)} \izo \ccO_{X_0} (d_1) \oplus \ldots \oplus \ccO_{X_0} (d_n). $$
Take $H\subset \P(V)$ to be an $n+1$ dimensional linear subspace general enough, so that it does not meet singular locus of $X$ (it is possible since the dimension 
of the singular locus in not greater than $n-3$, so $\dim H + \dim (X \backslash X_0) \leq n+1 + n-3 < 2n-1 = \dim \P(V)$). Define 
$$Y:= H \cap X = H \cap X_0 .$$
So Y is a projective curve and hence lemma \ref{o_wiazce_L} applies and therefore
$$
\xymatrix{
0 \ar@{->}[r]  &  T^*_{X_0}|_Y \ar@{->}[r]  &  \ccL|_Y \ar@{->}[r]^p                               &  \ccO_Y  \ar@{->}[r]                              &  0 \\
               &                            &  N_{X \slash \P(V)}(-2)|_Y \ar@{->}[u]_{\fii}^{\izo} &  \ccO_Y(d_i-2) \ar@{->}[l]_-{j_i} \ar@{.>}[u]_\xi        }
$$
Since $p \fii \ne 0$ (it is a composition of an isomorphism and an epimorphism) then there exists $i$ such that $\xi:=p \fii j_i \neq 0$, where $j_i$ is the inclusion 
of $i^{\textrm{th}}$ summand of the direct sum $\ccO_Y (d_1-2) \oplus \ldots \oplus \ccO_Y(d_n-2)$. Therefore $d_i -2 \leq 0$ and so either $d_i=1$
(which is the case, as I claim) or $d_i=2$ (which I need to exclude).

So suppose $d_i=2$. Then $\xi$ is an isomorphism and the composition $\fii j_i \xi^{-1}$ splits the exact sequence $ 0 \ra T^*_{X_0}|_Y \ra \ccL|_Y \ra \ccO_Y \ra 0$. 
Therefore $\ccL|_Y$ is the trivial extension, which contradicts the lemma \ref{o_wiazce_L}. Hence indeed $d_i=1$.
Now by the theorem \ref{o_rozm_zdegenerowanych} $X$ is a 
cone over a legendrian complete intersection of smaller dimension, again with singularities of codimension at least 2. So the induction on the dimension of $V$ can be applied.
\end{prf}

\begin{rem}
The assumption of the theorem is that $X$ is a normal complete intersection, while in the proof I only use that the
codimension of the singular locus is greater or equal to 2. In fact these are equivalent, due to 
\cite[\S17, thm 39]{matsumura}.
\end{rem}

\subsection{Examples}\label{examples}
At the very end of this section I present next few examples. The first three of them refer to the assumption of the 
theorem \ref{o_zupelnym_przecieciu}. 

\begin{ex}\label{four_lines}
Consider a union of four lines in $\P^3$ described by equations:
$$x_0 x_2 = 0 \ \textrm{ and } \ x_1 x_3 =0$$
All the lines are legendrian and so is their union. And they are described by two equations - but they are not irreducible.
\end{ex}

\begin{ex}\label{example_real}
Suppose that the ground field is $\R$. 
In $\P^5$ with coordinates \mbox{$(x_0:y_0:x_1:y_1:x_2:y_2)$} and symplectic form $\omega$ given by the matrix
$$\left(\begin{array}{rrrrrr}
 0 & 1 & 0 & 0 & 0 & 0\\
 -1& 0 & 0 & 0 & 0 & 0\\
 0 & 0 & 0 & 1 & 0 & 0\\
 0 & 0 & -1& 0 & 0 & 0\\
 0 & 0 & 0 & 0 & 0 & 1\\
 0 & 0 & 0 & 0 & -1& 0
\end{array} \right)$$
consider complete intersection of two quadrics and one cubic:
$$x_0^2 + y_0^2 - x_1^2 - y_1^2 = 0;$$ 
$$x_0^2 + y_0^2 - x_2^2 - y_2^2 = 0;$$ 
$$x_0x_1y_2 + x_0 y_1 x_2 + y_0 x_1 x_2 - y_0 y_1 y_2 = 0.$$ 
These equations look quite complicated, but if you think of $\P^5$ as a real projectivisation of a complex three dimensional vector space with coordinates $z_k = x_k + i y_k$,
where $i= \sqrt{-1}$ and $k=0,1,2$ they take much simpler form:
$$|z_0|^2 - |z_1|^2 =0;$$
$$|z_0|^2 - |z_2|^2 =0;$$
$$\operatorname{im} (z_0z_1z_2) =0. $$
Nevertheless, they define a legendrian variety and it is a smooth complete intersection.\footnote{This
 example was pointed out to me by Mark Gross - thanks!}
\end{ex}

\begin{ex}
Now take  take the same equations as in example \ref{example_real}, but defined over $\C$. The resulting surface
 is again legendrian and is a complete intersection, but it is not 
smooth  - it has singularities along six lines, so in codimension 1.
\end{ex}
\begin{prf}(sketch)
To verify legendrianity just use coordinates $x_k$ and $y_k$, then compute the conormal bundle and verify that $\omega'$ restricted to the conormal is 0. Use the theorem 
\ref{tangent_and_conormal}. If you skip to the section \ref{algebry_lie}, an easier way is to use the theorem \ref{Ideal_jest_podalgebraLie} and verify that the Lie bracket
of each pair of the above functions is 0. Conclude that the ideal of the variety is a Lie algebra using the remark \ref{ideal_generators}.

To compute singularities use coordinates $u_k: = x_k + i y_k$ and $v_k: = x_k - i y_k$ (which in real case could be written as $u_k:= z_k$ and $v_k:=\bar{z_k}$) so that 
the equations are:
$$u_0v_0 - u_1v_1 =0$$
$$u_0v_0 - u_2v_2 =0$$
$$u_0u_1u_2 - v_0v_1v_2=0$$
Verify, that the singular locus is an union of six lines:
$$u_0=v_0=u_1=v_2=0;$$ \nopagebreak
$$u_0=v_0=v_1=u_2=0;$$ \nopagebreak
$$u_1=v_1=u_0=v_2=0;$$ \nopagebreak
$$u_1=v_1=v_0=u_2=0;$$ \nopagebreak
$$u_2=v_2=u_0=v_1=0;$$ \nopagebreak
$$u_2=v_2=v_0=u_1=0;$$ \nopagebreak
\end{prf}

Before I present the last example, let me write an easy, though quite important lemma.
\begin{lem}\label{lemma_polynomial}
For a homogeneous polynomial $p$ in $n-1$ variables $y_1, \ldots y_{n-1}$ the following equality holds:
$$ y_1 \frac{\pd p}{\pd y_1} + \ldots + y_{n-1} \frac{\pd p}{\pd y_{n-1}}  =  \deg(p) \cdot p $$
\end{lem}
\begin{prf}
Compute for monomials and conclude for any homogeneous polynomial.
\end{prf}
\begin{ex}\label{general_example} (Compare with \cite[\S4.3]{landsbergmanivel04})
Let f be a homogeneous polynomial in $n-1$ variables $y_1, \ldots y_{n-1}$ of degree $k$. Denote by $X_f$ a subvariety in $\P^{2n-1}$ equal to the closure of the image of 
a map $\fii_f: \C^{n-1} \lra \P^{2n-1}$ defined by:
$$\fii_f(y):= \ \left(1:y_1: \ldots : y_{n-1}:(k-2)f(y): - \frac{\pd f}{\pd y_1}(y): \ldots: - \frac{\pd f}{\pd y_{n-1}}(y) \right)$$
where $y=(y_1, \ldots ,y_{n-1})$. Then $X_f$ is a legendrian subvariety for $\omega$ given by matrix 
$$J:= \left( \begin{array}{cc}
0&\Id_n\\
-\Id_n&0
\end{array} \right) .$$
\end{ex}
\begin{prf}
For short let me just write that $\fii_f(y):= \left(1:y:(k-2)f(y): - \ud f_y \right)$.
It suffices to check that for all $y \in \C^{n-1}$ the tangent space $T_y \hat{X_f}$ is Lagrangian. But $T_y \hat{X_f}$ is spanned by $n$ vectors:
$$u: =  \left(1,y,(k-2)f(y), - \ud f_y \right)  \textrm{ (direction tangent to the generator of the cone)}$$
$$v_i: = \frac{\pd \fii_f}{\pd y_1}(y) = \left(0,e_i,(k-2)\frac{\pd f}{\pd y_i}(y), - \frac{\pd (\ud f)}{\pd y_i}(y) \right) \textrm{ for } i =1, \ldots, k-1 $$
\indent \indent \indent \indent \indent \indent \indent \indent \indent  (directions tangent to the base of the cone)\\

\noindent where $e_i$ is an $i^{\textrm{th}}$ base vector of $\C^{n-1}$ and 
$$\frac{\pd (\ud f)}{\pd y_i}(y) = \left(\frac{\pd^2 f}{\pd y_1 \pd y_i}(y), \ldots,  \frac{\pd^2 f}{\pd y_{n-1} \pd y_i}(y) \right) .$$
So $\omega(v_i, v_j)$ is always 0 since $\frac{\pd^2 f}{\pd y_{i} \pd y_j}= \frac{\pd^2 f}{\pd y_{j} \pd y_i}$. Now compute $\omega(u, v_i)$:
$$\omega(u, v_i) =  (k-2)\frac{\pd f}{\pd y_i}(y) - y_1 \cdot \frac{\pd^2 f}{\pd y_1 \pd y_i}(y) -  \ldots -  y_{n-1} \frac{\pd^2 f}{\pd y_{n-1} \pd y_i}(y) + 
\frac{\pd f}{\pd y_i}(y) = $$
$$ =  (k-1)\frac{\pd f}{\pd y_i}(y) - \left( y_1 \cdot \frac{\pd}{\pd y_1} +  \ldots +  y_{n-1} \frac{\pd}{\pd y_{n-1}}\right) \left(\frac{\pd f}{\pd y_i}\right)(y) 
\stackrel{\textrm{(lemma \ref{lemma_polynomial})}}{=} 0$$
This proves that $X_f$ is legendrian. 
\end{prf}

Now let me write few words about the last example. For $\deg f = k =0, 1$ or $2$ it is not very interesting: in fact $X_f$ is just a linear 
subspace. So now assume $k=3$. Actually every smooth legendrian variety generated by quadrics is just $X_f$ for 
some choice of coordinates and for some $f$ of degree $3$. So let me now analyze some $X_f$ for $k=3$ and small $n$.\\

For $n=2$ the only (up to change of coordinates) non trivial polynomial of degree 3 is $f(y)=y^3$. In this case $X_f$ is the twisted cubic (see the example 
\ref{przyklad_skr_kubika}).\\

For $n=3$ I have few more polynomials: $f_1(y_1,y_2)=y_1^3$, $f_2(y_1,y_2)= y_1^2 y_2$ and $f_3(y_1,y_2) = y_1y_2(y_1+y_2)$.
 The first case is ''degenerated'' since $f_1$ 
does not depend on $y_2$. It is easy to see, that $X_{f_1}$ is contained in a hyperplane and it is just a cone over a twisted cubic. The second variety $X_{f_2}$ is
 isomorphic to product of $\P^1 \times \P^1$ embedded in $\P^5$ linearly on the first coordinate and quadratically on the second (see the example \ref{line_times_quadric}). 
The last case is not generated by quadrics, but there is quite a lot of them in its ideal. It is described by the following equations:
\begin{equation}
\left\{\begin{array}{l}
   x_0x_5 + x_1^2 + 2x_1x_2 =0;\\
   x_0x_4 + 2x_1x_2 + x_2^2  =0;\\
   3x_0x_3 + x_1x_4 + x_2x_5 =0;\\
   x_1x_4^2 - 2x_1x_4x_5 + 9x_2^2 x_3 - 5x_2x_4x_5 + 4x_2x_5^2=0; \\
   x_1x_3x_4 - 2x_1x_3x_5 + 2x_2x_3x_4 - x_2x_3x_5 - x_4^2x_5 + x_4x_5^2=0.\\
\end{array}\right.
\end{equation}
It has only one singular point at (0:0:0:1:0:0).\\

For $n=4$ things are getting much more complicated, since there are infinitely many projectively non-equivalent polynomials of degree 3. For some polynomials I can get very complicated
varieties (for example for $f=y_1^3 +y_2^3 +y_3 ^3$ $X_f$ is described by four quadrics, four quartics and one sixtic).\\

There are some questions concerning these examples: how to distinguish, for which $f$ the resulting variety $X_f$ 
is generated by quadrics? Are there any singular legendrian
varieties generated by quadrics? Conversely, if for $f$ of degree 3 the variety $X_f$ is smooth, 
does it imply that $X_f$ is generated by quadrics?


\section{Lie algebras}\label{algebry_lie}
In this section I describe a link between legendrian subvarieties in $\P^{2n-1}$ and some Lie algebras.

\begin{defin}
A linear space $\gotw$ together with a bilinear map $[\cdot , \cdot]: \gotw \times \gotw \ra \gotw$ is a \textbf{Lie algebra}, 
if the map (which is called the \textbf{Lie bracket}) is antisymmetric and it satisfies the \textbf{Jacobi identity}:
\begin{equation}
\big[w_1,[w_2,w_3]\big] + \big[w_2,[w_3,w_1]\big] + \big[w_3,[w_1,w_2]\big] = 0
\end{equation}
\end{defin}

Basic (and not only basic) properties of Lie algebras are in \cite[lectures 8-25]{fultonharris} or in \cite{humphreys}.

Now let me recall some classical examples of Lie algebras, which would be used later in this paper:

\begin{ex}
For a vector space $V$ ($\dim(V)=m$) define the algebra $\gotgl (V) = \gotgl_m$ as a vector space of endomorphisms of $V$ with the Lie bracket $[A,B]=A \circ B - B \circ A$.
\end{ex}
\begin{ex}
In $\gotgl (V)$ subalgebra of endomorphisms with trace equal to 0 is called $\gotsl (V)$ or $\gotsl_m$.
\end{ex}
\begin{ex}\label{przyklad_algebra_sp}
If in addition on $V$ ($\dim (V)= 2n$) there is defined a symplectic form $\omega$ then the symplectic Lie algebra $\gotsp (V) = \gotsp_{2n}$ is the space of 
linear endomorphisms $A: V \ra V$ such that
\begin{displaymath}
 \forall_{v,w \in V} \qquad \omega (Av,w) + \omega(v,Aw) =0.
\end{displaymath}
\end{ex}

If $\omega$ is a symplectic form on a vector space $V$ then it determines a structure of Lie algebra on 
the ring $\ccS=\F [V] = \bigoplus_i \Sym^i(V^*)$
of polynomial functions on $V$. To see that, recall the isomorphism $\phi$ defined in \eqref{izomorfizm_phi} and the induced form  $\omega':= (\phi^{-1})^* \omega$. 
This form $\omega'$ can be extended to a bilinear antisymmetric map:
$$\left[ \cdot , \cdot \right]: \ccS \times \ccS \ra \ccS$$
$$\left[ f , g \right] (x)= \omega' ( \ud f_x , \ud g_x)$$

Indeed it is an extension of $\omega'$: if $\alpha, \beta \in V^*=\Sym^1V^* \subset \ccS$, then
$[\alpha,\beta]=\omega'(\alpha, \beta)\in \C = \Sym^0V^* \subset \ccS$.\\

\begin{theo}(due to Poisson)
$\ccS$ together with the above map $[\cdot , \cdot]$ is a Lie algebra.
\end{theo}
\begin{prf}
Bilinearity of $[\cdot , \cdot]$ follows immediately from linearity of the derivation and from bilinearity of $\omega$. Antisymmetry is a consequence of the antisymmetry 
of $\omega$. It remains to verify the Jacobi identity.


Let $f,g,h \in \ccS$ be some polynomials. Then:\\
$$[ f \! \cdot \! g , \;  h ] (x)= \omega' \big( \ud (f_x \! \cdot \! g_x) , \;  \ud h_x \big) =  \omega' \big( f(x)\! \cdot \! \ud g_x + g(x) \! \cdot \! \ud f_x , \;  \ud h_x \big) =$$
$$ = f(x) \! \cdot \! \omega' ( \ud g_x , \; \ud h_x) + g(x) \! \cdot \! \omega' ( \ud f_x , \;  \ud h_x)= f(x) \! \cdot \! [ g , \;  h ] (x) + g(x) \! \cdot \! [ f , \;  h ] (x)$$\\
hence the bracket in $\ccS$ satisfies the Leibnitz formulae:
\begin{eqnarray}
 \left[f \! \cdot \! g , \;  h\right] = f\! \cdot \! \left[g , \;  h\right] +  g\! \cdot \! \left[f , \;  h\right] \nonumber \\
 \left[f, \; g\! \cdot \! h\right] = g\! \cdot \! \left[f, \;
h\right] + h\! \cdot \! \left[f, \; g\right].
\label{rown_Leibnitza}
\end{eqnarray}

Now let me define a trilinear map $\alpha: \ccS \times \ccS \times \ccS \ra \ccS$
\begin{displaymath}
\alpha(f, \; g, \; h):= \Big[f, \; [g, \; h]\Big] + \Big[g, \; [h,
\; f]\Big] + \Big[h, \; [f, \; g]\Big].
\end{displaymath}
The Jacobi identity is verified if I prove that $\alpha \equiv 0$. So let me use (\ref{rown_Leibnitza}) for the bracket to show similar property for $\alpha$ and then 
I will have to verify the identity only for the multiplicative generators of $\ccS$ (i.e. linear functions).\\

If I substitute $f=f_1 \! \cdot \! f_2$ into the defining equation of $\alpha$ and then apply (\ref{rown_Leibnitza}) five times, I will get ten summands of which there are
two pairs of opposite terms. So there will remain six entries, three of  which can be again composed into 
\mbox{$f_1 \! \cdot \! \alpha (f_2, \;  g, \;  h)$} and the other three into 
\mbox{$f_2 \! \cdot \! \alpha (f_1, \;  g, \;  h)$}. So:
$$\alpha(f_1 \! \cdot \! f_2, \; g, \; h)= 
 f_1 \! \cdot \! \alpha (f_2, \;  g, \;  h) + f_2 \! \cdot \! \alpha (f_1, \;  g, \;  h).$$

And in the same way one can prove that:
$$\alpha(f, \;  g_1 \! \cdot \! g_2, \;  h) =  g_1 \! \cdot \! \alpha (f, \;  g_2, \;  h) + g_2 \! \cdot \! \alpha (f, \;  g_1, \; h) $$
$$\alpha(f, \;  g, \;  h_1 \! \cdot \! h_2) =  h_1 \! \cdot \! \alpha (f, \;  g, \;  h_2) + h_2 \! \cdot \! \alpha (f, \;  g, \; h_1) $$

Hence if only $\alpha(f_1, g, h) = \alpha(f_2, g, h) = 0$, then $\alpha(f_1f_2,g,h)=0$ as well (and similarly for $g_1$, $g_2$ and $h_1$, $h_2$).
It remains to prove that $\alpha = 0$ just for multiplicative generators of the ring $\ccS$, so for linear polynomials. But if $f,g,h$ are linear, then \mbox{$\large[f,[g,h]\large]=0$},
since the second differential of a linear function is 0. Hence \mbox{$\alpha(f,g,h)=0$.}
\end{prf}

Although the fact above itself seems to be very interesting on it's own, it is not clear yet what does it has to do
 with the legendrian subvarieties of $\P(V)$. Here follows 
the main theorem of this section which shows that indeed there is a correspondence between some Lie subalgebras of $\ccS$ and legendrian subvarieties of $\P(V)$ 

\begin{theo}\label{Ideal_jest_podalgebraLie}
Suppose that $X \subset \P(V)$ and that $\I=\I(X) \ideal \ccS$ is the ideal describing $X$. Then the following conditions are equivalent:
\begin{itemize}
\item[(i)] 
$X$ is legendrian,
\item[(ii)]
$\I$ is Lie subalgebra of $\ccS$ and each irreducible component of $X$ is $n-1$ dimensional.
\end{itemize}
\end{theo}

\begin{prf}
This theorem is known in the theory of D-modules (see 
\cite[chapter 11,.prop. 2.4]{coutinho})\footnote{This
 remark is due to Mircea Mustata.}.
Let me rewrite the proof so that it is more suitable for the context of this paper.\\

By the theorem \ref{tangent_and_conormal} I know that (i) is equivalent to (i'):
\begin{itemize}
\item[(i')] $\dim \hat{X} =n$ and $\omega'|_{N^*_{\hat{X}_0\slash V}} \equiv 0$.
\end{itemize}

So it suffices to show that $\omega'|_{N^*_{\hat{X}_0\slash V}} \equiv 0$ if and only if $\I$ is a Lie subalgebra in $\ccS$.\\

Suppose that $x \in \hat{X}_0$ is any point and that $f,g \in \I$ are any polynomials vanishing on $\hat{X}$. Then surely $\ud f_x|_{T \hat{X}_0} \equiv 0$, so
$\ud f_x \in N^*_{\hat{X}_0 \slash V}$ and similarly $\ud g_x \in N^*_{\hat{X}_0 \slash V}$.\\

If I suppose $\omega'|_{N^*_{\hat{X}_0\slash V}} \equiv 0$ then
\begin{displaymath}
[f,g](x) = \omega'(\ud f_x, \ud g_x) = 0,
\end{displaymath}
i.e. $[f,g]|_{\hat{X}_0} = 0$, so $[f,g]|_{\hat{X}} = 0$ (because the closure of $\hat{X}_0$ is exactly $\hat{X}$) and so
  $[f,g] \in \I$.\\

Conversely, if $\I$ is a Lie subalgebra then
\begin{displaymath}
\omega'(\ud f_x, \ud g_x )= [f,g](x) = 0.
\end{displaymath}
Since the map
$$\begin{array}{ccl}
  \I & \lra & N^*_{x,\hat{X}_0\slash V} \\
  f & \longmapsto & \ud f_x
\end{array}$$
is an epimorphism of vector spaces for each $x \in \hat{X}_0$ then $\omega'|_{N^*_{\hat{X}_0\slash V}} \equiv 0$.
\nopagebreak
\end{prf}

\begin{rem}\label{ideal_generators}
Suppose $\I \ideal \ccS$ and $f,g \in \I$, $s \in \ccS$. If $[f,g] \in \I$ then also $[sf,g] \in \I$. In particular, if $\I=(f_1, f_2, \ldots, f_k)$ (i.e. $\I$ is generated as a ring ideal by $f_i$'s) and for each 
$i, j \in \{1 \ldots k\}$ I have $[f_i,f_j] \in \I$, then $\I$ is a Lie subalgebra.
\end{rem}

\begin{prf}
It follows directly from the equation \eqref{rown_Leibnitza}.
\end{prf}

The ring $\ccS$ admits canonical gradation $\ccS=\bigoplus_i \Sym^i V^*$. Suppose $f \in \Sym^i V^*$ and $g \in \Sym^j V^*$ (so $f$ and $g$ are homogeneous of degree 
respectively $i$ and $j$). Since each of the derivations decreases the degree by one, it follows that
\begin{equation}\label{formula_na_gradacje}
[f,g]\in \Sym^{i+j-2} V^*.
\end{equation}

\begin{cor}\label{Sym2_jest_podalg}\hfill\nopagebreak
\begin{itemize}
\item[(a)] $\Sym^2 V^*$ is a finite dimensional Lie subalgebra in $\ccS$.
\item[(b)] $\Sym^2 V^*$ is isomorphic (as a Lie algebra) with $\gotsp(V)$ (see example \ref{przyklad_algebra_sp}).
\item[(c)] If $X$ is legendrian subvariety in $\P(V)$ and $\I = \I(X)$, then $\I_2:= \I \cap \Sym^2 V^*$ is a Lie subalgebra of $\Sym^2 V^*$.
\end{itemize}
\end{cor}

\begin{prf}
First notice that (a) follows from formulae \eqref{formula_na_gradacje} with $i=j=2$. Next (c) follows immediately from (a) and from theorem \ref{Ideal_jest_podalgebraLie}.
It remains to prove (b).\\

Now let $\{e_1, \ldots, e_n, \ e_{n+1}, \ldots, e_{2n}\}$ be a standard (due to $\omega$) basis of $V$, i.e. in this basis $\omega$ corresponds to the matrix $J$ where
$$J=\begin{pmatrix}
  0 & \Id_n \\
  -\Id_n & 0
\end{pmatrix}.$$
($\Id_n$ is the $n \times n$ identity matrix)

A choice of basis identifies the vector space $\Sym^2 V^*$ with the space of symmetric matrices $2n \times 2n$. So a quadratic polynomial $f \in \Sym^2 V^*$ is identified 
with such a matrix $A$, that
\begin{displaymath}
 f(x) \: = \: x^T A x.
\end{displaymath}
Let me compute what is the matrix of the Lie bracket of two polynomials. Suppose $f(x)=x^T A x$, $g(x)=x^T B x$. An easy calculation proves that the matrix of $\omega'$
in the dual basis is again $J$. So:
$$\ud f_x = 2x^T A \ \textrm{ and } \ \ud g_x = 2x^T B;$$
$$[f,g](x) \; = \; \omega'(\ud f_x, \ud g_x) \; = \; (2x^T A) J
(2Bx) \; = \; x^T (4AJB) x \; =$$
$$= \; x^T \Big(2\big(AJB \ + \ (AJB)^T \big)\Big)x \;
= \; x^T \big(2(AJB \ - \ BJA) \big)x$$
 Hence the Lie bracket in the space of symmetric matrices is defined by:
\begin{equation}\label{nawias_Lie_w_Sym2}
[A,B] \; = \; 2(AJB \ - \ BJA)
\end{equation}
 Now I can easily define an isomorphism of $\Sym^2 V^*$ and $\gotsp(V)$ being just the multiplication by $2J$:
\begin{equation}\label{isomorphism_Sym2_sp}
\begin{array}{rrcl}
  \rho :&\Sym^2 V^*&\lra       &\gotsp(V) \\
        &f\izo A &\longmapsto& 2JA
\end{array}
\end{equation}

Again simple calculation shows that $\rho$ is indeed a linear isomorphism. 
Moreover th equation \eqref{nawias_Lie_w_Sym2} proves that $\rho$ preserves the Lie bracket, so it is an isomorphism of Lie algebras.
\end{prf}

In fact, one can easily prove, that the isomorphism $\rho$ does not depend on the choice of the symplectic basis 
$\{e_1 \ldots e_{2n}\}$.


\section{Symplectic group} \label{rozdzial_rozm_gen_kwadryki}

Theorem \ref{o_rozm_zdegenerowanych} explained what happens, if there is a linear polynomial in the ideal $\I(X)$.
In this subsection I study the quadratic part $\I_2$ if the ideal.
It gives rise a group of projective automorphisms of X.
 
\begin{defin}
A \textbf{symplectic group} of $V$ (denoted $\Sp$) is the group of linear automorphisms of $V$ preserving the symplectic form $\omega$, i.e. $\psi \in \Sp$ 
if and only if $\omega(\psi(v), \psi(w)) = \omega(v,w)$ or equivalently $\psi^T J \psi = J$ where $J$ is the matrix of $\omega$.
\end{defin}

Finite dimension of the Lie algebra $\Sym^2 V^* \izo \gotsp(V)$  makes possible passing to Lie groups (see \cite[chapters II-III]{humphreys2} and 
\cite[section 8.3]{fultonharris}). To be precise, I have the exponential map $\exp : \gotsp(V) \lra \Sp$ which in particular takes Lie subalgebra $\gotg$ into
unique connected subgroup of $\Sp$ such that its tangent space is exactly $\gotg$.\\

At the begining compare two representations of $\Sp$. The first is the right $\Sp$ action on $\gotsp(V)$:
\begin{displaymath}
\begin{array}{rrcl}
 \Ad:&\Sp       & \lra      & \Gl (\gotsp(V))\\
     &\psi      &\longmapsto& \psi^{-1} \circ \ \cdot \ \circ \psi \textrm{, i.e.}\\
     &A\star\psi& =         & \psi^{-1} \circ A \circ \psi,
\end{array}
\end{displaymath}
where $\circ$ is just the composition of linear maps and $\star$ is the group action. This is the adjoint action.\\

The other representation is the natural representation, i.e. ~a left action of $\Sp$ on $V$: an element $\psi$ of the group takes a vector $v$ simply to the vector $\psi(v)$.
Let me denote this representation by $\tau : \Sp \mono \Gl(V)$.\\

The representation $\tau$ inducts a right action $\bullet$ of the group $\Sp$ on the ring of polynomials $\ccS$: $(f \bullet \psi) (x) = f(\psi(x))$ or simply 
$f \bullet \psi = f \circ \psi$ $\circ$ is just map composition). This action preserves gradation so it restricts to a right action $\tau^i$ on $\Sym^i V^*$.

\begin{theo} \label{reprezentacje_sa_rowne}
Representations $\Ad: \Sp \ra \Gl(\gotsp(V))$ and $$\tau^2: \Sp \lra \Gl(\Sym^2(V)) \stackrel{\rho}{\izo} \Gl(\gotsp(V))$$
(where $\rho$ is defined as in the equation \eqref{isomorphism_Sym2_sp}) are equal.
\end{theo}

\begin{prf}
Choose a basis of $V$ standard due to $\omega$. Let $\psi \in \Sp$ be an arbitrary element. I will check what automorphism of $\gotsp(V)$ is $\tau^2(\psi)$.
So let $f \in \Sym^2(V^*)$ be an arbitrary quadratic polynomial and suppose $f(x) = x^T A x$ for a symmetric matrix $A$. Then 
$(f \bullet \psi)(x)=f(\psi(x))= x^T \psi^T A \psi x$
 and hence $f \bullet \psi$ corresponds to the matrix $\psi^T A \psi$. So under the isomorphism \mbox{$\rho: \Sym^2 V^* \ra \gotsp(V)$} 
(see \eqref{isomorphism_Sym2_sp}) $f \bullet \psi$ goes to $2J \psi^T A \psi$.
But since $\psi^T J \psi = J$  and $J^2 = -\Id_{2n}$, then $J \psi^T = - J^{-1} \psi^T = - \psi^{-1} J^{-1} = \psi^{-1} J$, so
$$\rho \circ \tau^2(\psi)(f) = \rho(f \bullet \psi) = 2J \psi^T A \psi = \psi^{-1} 2 J A \psi = \psi^{-1} \rho(f) \psi.$$
So $\rho \circ \tau^2 = \Ad$.
\end{prf}

\begin{theo} \label{o_rozniczce_tau_i}
For any $i\in \N$ the derivation of the following homomorphism of Lie groups:
\begin{displaymath}
\begin{array}{rrcl}
\tau^i:& \Sp& \lra & \Gl(\Sym^i V^*)\\
&\psi & \longmapsto& (f \mapsto f\bullet \psi)
\end{array}
\end{displaymath}
is the homomorphism of Lie algebras:
\begin{displaymath}
\begin{array}{rrcl}
\uD \tau^i:& \Sym^2V^*& \lra & \gotgl(\Sym^i V^*)\\
&f & \longmapsto& [\cdot, f]
\end{array}
\end{displaymath}
\end{theo}

\begin{prf}
The set of maps $\uD \tau^i$ satisfies the following conditions:
\begin{itemize}
\item $\uD \tau^2= \uD \Ad = \ad $ (this follows from the theorem \ref{reprezentacje_sa_rowne})
\item (Leibnitz rule) $\uD \tau^{i+j}(\cdot)(fg) = \uD \tau^i(\cdot)(f) \cdot g + f \cdot \uD \tau^j(\cdot)(g) $ for each $f\in \Sym^i V^*$, $g\in \Sym^j V^*$
\end{itemize}
First I prove that $\uD \tau^1= \ad$. Indeed, for any $\alpha \in V^*$ I have:
\begin{displaymath}
2\alpha \ \uD \tau^1(\cdot)(\alpha)= \uD \tau^2(\cdot)(\alpha^2) = [\cdot,\alpha^2]= 2\alpha [\cdot, \alpha]. 
\end{displaymath}
So $\uD \tau^1(\cdot)(\alpha)= [\cdot,\alpha]$. Now applying the Leibnitz rule for both derivation and Lie bracket proves the theorem.
\end{prf}

\begin{cor}\label{G_zachowuje_g}
Let $\I$ be any Lie subalgebra of $\ccS$ and let $\I_2 = \I \cap \Sym^2V^*$. Suppose that $G < \Sp$ is the subgroup corresponding to $\I_2 \subset \Sym^2V^* \izo \gotsp(V)$. Then the action $\tau^*$
restricted to $G$ preserves $\I$.
\end{cor}

\begin{prf}
Since $\I_2$ is a subalgebra of $\I$ surely its adjoint action preserves $\I$. But I have shown in the theorem \ref{o_rozniczce_tau_i} that $\uD \tau^*$ is equal to 
the adjoint action. So $\tau^*$ preserves $\I$ as well.
\end{prf}

\begin{cor}\label{G_zachowuje_X}
Suppose $X \subset_l \P(V)$  and $\I=\I(X)$ is the ideal describing $X$. Let $\I_2= \I \cap
\Sym^2 V^*$ and $G < \Sp$ be the corresponding subgroup (as in the corollary \ref{G_zachowuje_g}). Now $G$ acts on $V$ (by natural representation $\tau|_G$)
and the action preserves $\hat{X}$ (so also $G$ acts on $\P(V)$ and it preserves $X$).
\end{cor}

\begin{prf}
By the theorem \ref{Ideal_jest_podalgebraLie} $\I$ is a Lie subalgebra in $\ccS$ so the group $G$ is well defined in the statement of the corollary.
It follows from the corollary \ref{G_zachowuje_g} that $\tau^*|_G$ preserves $\I$. So $\tau|_G$ preserves $\hat{X}$. Since $G$ acts by linear automorphisms it acts on 
$\P(V)$ and this action preserves $X$.
\end{prf}

\begin{lem}\label{G_is_closed}
Let $X$, $\I$, $\I_2$ and $G$ be as in corollary \ref{G_zachowuje_X}. Then $G$ is a maximal connected subgroup in $\Sp$ 
preserving $X$. In particular $G$ is closed and so it is a
Lie group.
\end{lem}

\begin{prf}
Suppose $H$ is a connected subgroup of $\Sp$ containing $G$ and preserving $X$. I have to show that $H=G$.

Let $\J \subset \Sym^2V^*$ be the Lie algebra corresponding to $H$. Then $\tau^*|_H$ preserves $\I$ and so by theorem \ref{o_rozniczce_tau_i}:
$$[\J,\I] \subset \I$$
Fix a homogeneous polynomial $f\in \J$, choose an arbitrary $x\in \hat{X}_0$ and take any $g \in \I$. Since $[f,g] \in \I$:
$$0=[f,g](x)=\omega'(\ud f_x, \ud g_x).$$
This happens for all $g\in \I$ so $\ud f_x$ is perpendicular (with respect to $\omega'$) to $N^*_{x,\hat{X}_0 \slash V}$ so $\ud f \in N^*_{x,\hat{X}_0 \slash V}$ and 
$\ud f_x |_{T_x\hat{X}_0} = 0$. Since $x$ was chosen arbitrary, it follows that $\ud f|_{T \hat{X}_0} \equiv 0$. So $f$ is constant on $\hat{X}_0$ 
(and on $\hat{X}$ as well). But $f$ is homogeneous and hence simply $f \in \I$. Again $f$ was chosen arbitrary, so $\J=\I$ and $H=G$. This proves maximality of $G$.

If $\bar{G}$ is the closure of $G$, then it preserves $X$. So $\bar{G}=G$ i.e. $G$ is closed. A closed subgroup of a Lie group is again a Lie group and so is G.
\end{prf}

One more is interesting to notice: 
whenever $X$ is irreducible and non-degenerate, the action of $G$ is (almost) faithful.

\begin{prop}\label{kernel_of_action}
If $X$ is non-degenerate and irreducible then the subgroup of $G$ acting trivially on $X$ is either trivial or 
$\Z\slash 2\Z$, hence discrete.
\end{prop}

\begin{prf}
Assume $g \in \Sp$ acts trivially on whole of $X$. The locus of points in $\P(V)$ fixed by $g$ is just the (disjoint)
 union of the eigenspaces of $g$. Since $X$ is irreducible, it is contained in one of the eigenspaces, which must be the
whole of $\P(V)$, because $X$ is nondegenerate. So in fact $g= \lambda \Id$ for some $\lambda \in \C^*$. Since
$g^T J g = J$, it follows that $\lambda^2=1$, hence $g=\pm \Id$.
\end{prf}

Now let me explore the action of $G$ on $X$ more precisely. If $X$ is a general enough little can be deduced from this
 action, since $\I_2$ (and so $G$) might be very
small (for example trivial). But if I restrict to legendrian varieties which are generated by quadrics 
(i.e. $\I_2$ generates whole the ideal $\I$) there are several
interesting results, for example that there exist a huge open orbit of this action: it is $X_0$ itself! 
In particular, if in addition $X$ is smooth it follows
that $X$ is homogeneous and a full  classification might be given.

\begin{theo} \label{przeciecie_X'_i_X_0_jest_puste}
Assume that $X \subset_l \P(V)$ is irreducible and $\I(X)$ is generated by quadrics 
and suppose that $G$ is the group acting on $X$ as defined in corollary \ref{G_zachowuje_X}.
Also let $Y\varsubsetneq X$ be a closed subvariety which is invariant under action of $G$. Then $Y\cap X_0 = \emptyset$.
\end{theo}

\begin{prf}
Let $\J = \I(Y) \ideal \ccS$. Since G preserves $Y$, it preserves $\J$ as well, so by the theorem \ref{o_rozniczce_tau_i}:
\begin{displaymath}
[\J,\I_2] \subset \J
\end{displaymath}
and from the remark \ref{ideal_generators} also $[\J,\I] \subset \J$. 

Now I will simulate the proof of \mbox{the lemma \ref{G_is_closed}:}
suppose there exists a point $x \in \hat{X_0}\cap \hat{Y}$ and fix an arbitrary $f\in \J$. Then for all $g \in \I$
\begin{displaymath}
0 \; = \; [f,g](x) \; = \; \omega'(\ud f_x, \ud g_x)
\end{displaymath}
so $\ud f \in (N^*_{x,\hat{X}_0 \slash V})^{\perp_{\omega'}} = N^*_{x,\hat{X}_0 \slash V}$.
It holds for every $f \in \J$, hence the dimension of conormal space to $\hat{Y}$ at point $x$ is $n$.
But the set $\hat{X_0}\cap \hat{Y}$ is open in $\hat{Y}$. So $\dim \hat{Y} = 2n - n = n = \dim \hat{X}$. Since $Y$ is closed and $X$ is irreducible, $Y=X$ although 
I have assumed that $Y$ is proper subvariety of $X$. So this contradicts the assumption that there exists a point $x \in \hat{X}_0 \cap \hat{Y}$.
So $\hat{X}_0 \cap \hat{Y} = \emptyset$ and $X_0 \cap Y = \emptyset$ as well.
\end{prf}

\begin{cor} \label{orbit_X_0}
If $X\subset_l \P^{2n-1}$ is irreducible and generated by quadrics then $X_0$ is an orbit of the action of $G$.
\end{cor}
 
\begin{prf}
Group $G$ from the corollary \ref{G_zachowuje_X} acts on $X$. Suppose $\ccO$ is an orbit of any smooth point of $X$. The closure $\bar{\ccO}$ of the orbit is invariant under
the action of $G$, so from the theorem \ref{przeciecie_X'_i_X_0_jest_puste} it follows $\bar{\ccO}=X$. 
Also \mbox{$\bar{\ccO}\backslash \ccO$} is invariant so again from the theorem \ref{przeciecie_X'_i_X_0_jest_puste} 
\mbox{$(\bar{\ccO}\backslash \ccO) \cap X_0 = \emptyset$}, i.e. $\ccO = X_0$. 
\end{prf}

\begin{theo} \label{rozmaitosci_jednorodne}
If $X\subset_l \P^{2n-1}$ is smooth irreducible and generated by quadrics then it is a homogeneous space.
\end{theo}
\noprf

\begin{rem}\label{real_remark}
Most of the things I have done so far, works also over the field which is not algebraically closed (but still of 
characteristic 0), for example for the real numbers. The only exception is the subsection \ref{section_proj_geometry} 
in particular the lemma \ref{o_wiazce_L} and the theorem \ref{o_zupelnym_przecieciu}. But in the section 
\ref{quadrics} I will heavily use the the assumption of the ground field being $\C$. Nevertheless, I believe with some 
(but only a little) effort one could follow very similar argument for the reals. Anyway, at the moment I cannot see any
reasonable application for such theory over reals and for that reason I will not complicate the text.
\end{rem}



\section{More examples - subadjoint varieties}  \label{rozdzial_przyklad}

\subsection{Digression on contact manifolds}

\begin{defin}
A complex projective manifold $Y^{2n+1}$ (of dimension $2n+1$) with fixed a rank $2n$ subbundle of the tangent bundle 
$F^{2n}\subset T X$ is a \textbf{contact manifold} if the  form $F \times F \ra TX \slash F =:L$ induced by a Lie bracket 
is nowhere degenerate. In particular this means that there is a symplectic form on every fibre of $F$
(for some equivalent descriptions see for example \cite{wisniewski}).
\end{defin}

The contact manifolds come up from the classification of quaterionic-k\"{a}hler manifolds - the twistor space of 
such turns out to be a complex contact manifold. The detailed description and original references can be found in 
a review on the subject \cite{beauville}.\\

The four authors \cite{wisniewski} proved the following theorem:
\begin{theo}
If $Y$ is a complex contact manifold then one of the following holds:
\begin{itemize}
\item $Y$ is a Fano variety with second Betti number $b_2=1$ or 
\item $Y$ is a (Grothendieck) projectivisation of the tangent bundle to some projective manifold $Y'$ or
\item the canonical divisor $K_Y$ is numerically effective.
\end{itemize}
\end{theo}
The last case has been excluded by \cite{demailly}. The interesting case is the first one: 

\begin{con} \label{hipoteza_contact}
If $Y$ is a Fano complex contact manifold with $b_2 = 1$ then $Y$ is a homogeneous variety which is the unique 
closed orbit of the adjoint action of some simple Lie group $G$ on $\P(\gotg)$ 
(where the $\gotg$ is the Lie algebra tangent to $G$).
\end{con}

These varieties are called minimal nilpotent orbits 
(see \cite{beauvillefano}, \cite{beauville}, \cite{wisniewski}) or adjoint varieties (see \cite{landsbergmanivel04}). 
They were put in the table \ref{table_contacts}.

The attacks on the conjecture invoked  the legendrian varieties (see: \cite{contact}, \cite{kebekus1}, \cite{kebekus2}).
They are related to the notion of contact line:
\begin{defin}
For a Fano contact manifold $Y$ with $b_2=1$, a rational curve  \mbox{$C\subset Y$} is \textbf{contact line}, if
 $L\cdot C=1$ (where  $L=TY \slash F$).
\end{defin}

Now, if $y \in Y$ is a general point of a Fano contact $Y$  then the tangent directions to the contact lines 
passing through $y$ generate a smooth legendrian subvariety $X_y$ in  $\P(F_y)$.
\footnote{There is only one exception,
which in a sense is degenerate: on $\P^{2n+1}$ the line bundle $L$ is isomorphic to $\ccO(2)$ and therefore 
there are no contact lines on $\P^{2n+1}$ and the 'legendrian variety' is empty.}

\begin{table}[hbt]
\label{table_contacts}
{\small\centering
\begin{tabular}{|p{0.1\textwidth}|p{0.1\textwidth}|p{0.21\textwidth}|p{0.17\textwidth}|p{0.25\textwidth}|}
\hline
  Lie group  & type& contact  manifold $Y^{2n+1}$     & legendrian variety $X^{n-1}$ & remarks \\
\hline
&&&&\\
$\operatorname{SL}_{n+2}$  & $ A_{n+1}  $              & $\P(T\P^{n+1})$      & $\P^{n-1} \sqcup \P^{n-1}$ \mbox{$\subset \P^{2n-1}$}
& $b_2(Y) =2 $    \\
\hline
$\operatorname{Sp}_{2n+2}$ & $ C_{n+1}$              & $\P^{2n+1}$          &
 $\emptyset\subset \P^{2n-1}$                & $Y$ does not have any lines\\
\hline
$\operatorname{SO}_{n+4}$  & $ B_{\frac{n+3}{2}}$  or $D_{\frac{n+4}{2}}$  & $Gr_O(2,n+4)$        & 
$\P^1 \times Q^{n-2}$ \mbox{$ \subset \P^{2n-1}$}   & $Y$=Grassmannian  of projective lines on a quadric $Q^{n+2}$\\
\hline
                           & $ G_2$                  &Grassmannian of special lines on $Q^5$ & $\P^1 \subset \P^3$  & $X$=twisted cubic\\
\hline
                           & $ F_4$                  & an $F_4$ variety   & $Gr_L(3,6)$  \mbox{$ \subset \P^{13}$}    &   \\
\hline
                           & $ E_6$                  & an $E_6$ variety   & $Gr(3,6) \subset \P^{19}$        &   \\
\hline
                           & $ E_7$                  & an $E_7$ variety   & $\mathbb{S}_6 \subset \P^{31}$  &$X$=spinor variety\\
\hline
                           & $ E_8$                  & an $E_8$ variety   & an $E_7$ variety $\subset \P^{55}$&   \\
\hline
\end{tabular}
\caption{Simple Lie groups and corresponding minimal nilpotent orbits ($Y$) together with theirs varieties of
directions tangent to lines ($X$).}}
\end{table}

The natural question arises: which among the legendrian varieties arise in this way?
For the adjoint varieties (which are the only known examples) one gets the list of homogeneous legendrian varieties, 
called \textbf{the subadjoint varieties} (see \cite{landsbergmanivel}, \cite{mukai}). They are expected to be the only 
homogeneous legendrian varieties (a partial proof can be found in \cite{landsbergmanivel04}) and it known that
they are the only symmetric legendrian varieties. Another characterisation comes out in this paper: they are the only 
smooth irreducible legendrian varieties, whose ideal is generated by quadratic polynomials (see theorem 
\ref{final_classification}). Jaros\l{}aw Wi\'sniewski hopes to prove that every legendrian variety arising from the 
Fano contact manifolds with $b_2=1$ is in fact generated by quadrics. It is hoped to be a way to prove the contact conjecture.

In this section I am going to present all the subadjoint varieties.

\subsection{Line times quadric} \label{line_times_quadric}
For the groups $\mathbf{SO}_{n+4}$ 
(i.e of type $B_{\cdot}$ and $D_{\cdot}$) one gets
 the adjoint variety $G_o (2,n+4)$ (i.e. the Grassmannian of lines on a quadric hypersurface).
 Its legendrian variety is a product of a line and a quadric hypersurface.

\begin{figure}[htb]
\centering
\includegraphics[width=0.9\textwidth]{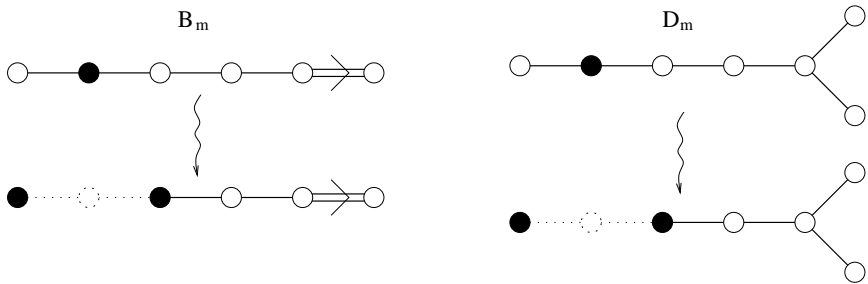}
\caption{The two homogeneous varieties (the minimal nilpotent orbit for $B_m$ or $D_m$ (above) and its variety of lines
 (below)) represented on the Dynkin diagrams. The black dots denote the simple roots that are not in the corresponding 
parabolic subgroup.}
 \label{pic_adjoint_Bn_and_Dn}
\end{figure}

For any \mbox{$n \geq 3$} consider a quadric hypersurface $Q^{n-2} \subset \P^{n-1}$ described by the equation \mbox{$y^Ty=0$.}\footnote{
It is convenient to chose such a quadric in order to get a uniform 
description of both even and odd dimensional cases. Yet in the context of 
representation theory and the subject of section \ref{quadrics} a different
choice should be proposed.}
Take  the Segre embedding $\fii: \P^1 \times \P^{n-1} \mono \P^{2n-1}$
\begin{displaymath}
\fii((\lam:\m),(y_0:\ldots:y_{n-1}))= (\lam y_0:\ldots: \lam y_{n-1}: \m y_0: \ldots: \m y_{n-1} ).
\end{displaymath}
Then $X:=\fii(\P^1 \times Q^{n-2})$  is a legendrian subvariety in $\P^{2n-1}$.\\

This can be verified using the theorem \ref{Ideal_jest_podalgebraLie}. Let me start with computing the ideal of $X$. Let 
\begin{itemize}
\item $f_{ij}:= x_i x_{n+j} - x_{j} x_{n+i}$ for $i,j = 0,1,\ldots ,n-1$
\item $g_{+}:= \half (x_n^2 + x_{n+1}^2 \ldots + x_{2n-1}^2)$
\item $g_{-}:= -\half (x_0^2 + x_1^2 \ldots + x_{n-1}^2)$
\item $h := x_0 x_n + x_1 x_{n+1} + \ldots + x_{n-1} x_{2n-1}$
\end{itemize}

Functions $f_{ij}$ generate the ideal of whole $\P^1 \times \P^{n-1}$ and the ideal of $X$  is $\I=(\{f_{ij}\}_{i,j=0,\ldots n-1}, \: g_+, g_-, h)$. Since 
$\dim X = n-1$, it remains to verify if $\I$ is a Lie subalgebra. 
 Notice that $f_{ij}= - f_{ji}$. Then  compute the derivations of $f_{ij}$ and 
$g_+$, $g_-$, $h$:
\begin{itemize}
\item $\ud f_{ij} = x_{n+j} \ud x_i - x_{n+i} \ud x_j - x_j \ud x_{n+i} + x_i \ud x_{n+j}$
\item $\ud g_{+} = (0, 0, \ldots, 0, \ x_n, x_{n+1}, \ldots , x_{2n-1})$
\item $\ud g_{-} = -(x_0, x_1,\ldots , x_{n-1}, \  0,0, \ldots, 0 )$
\item $\ud h = (x_n, x_{n+1},\ldots , x_{2n-1}, \ x_0, x_1 \ldots, x_{n-1} )$
\end{itemize}
Now here are the Lie brackets of the generators of $\I$:
(form $\omega'$ is standard, like usually):
\begin{itemize}
  \item[(i)] $[f_{ij},f_{jk}]= - x_{n+i} x_{k} + x_{i} x_{n+k} =
  f_{ik}$ for $i\ne j \ne k$
  \item[(ii)] $[f_{ij},f_{kl}]= 0 $ if $i,j,k,l$ are different numbers $\in \{ 0,1, \ldots n-1\}$
  \item[(iii)] $[f_{ij},g_+]= x_{n+j} x_{n+i} - x_{n+i} x_{n+j} = 0$
  \item[(iv)] $[f_{ij},g_-]= -x_j x_i + x_i x_j = 0$
  \item[(v)] $[f_{ij},h]= x_{n+j} x_{i} - x_{n+i} x_{j} + x_j x_{n+i} - x_i x_{n+j} = 0$
  \item[(vi)] $[g_+, g_-]= x_n x_0 + x_{n+1} x_1 + \ldots + x_{2n-1}x_{n-1} = h$
  \item[(vii)] $[h, g_-]=  x_0^2 +x_1^2 + \ldots + x_{n-1}^2 = -2 g_-$
  \item[(viii)] $[h, g_+]=  x_n^2 + x_{n+1}^2 + \ldots + x_{2n-1}^2 = 2 g_+$
\end{itemize}

So indeed $\I$ is a Lie subalgebra and $X$ a legendrian subvariety.\\

The next thing to do is to understand what kind of algebra is $\I_2$. As a linear space, it is spanned by $ \Big\{ \{ f_{ij} \}_{i,j=0,\ldots n-1}, \ g_+, g_-, h \Big\} $,
where $f_{ij}=-f_{ji}$. Equations (i)-(viii) show, that $\I_2$ is a direct sum of the subalgebra spanned by $\{f_{ij}\}$ and the subalgebra spanned by $\{g_+, g_-, h\}$.
It is clear that the second one is isomorphic to $\gotsl_2$. The other one is actually isomorphic to $\gotso_n$, 
which can be viewed as an algebra of skew-symmetric 
$n\times n$ matrices. But then if I define $f_{ij}$ to correspond to the elementary skew matrix with 1
 at the position $(i,j)$, $-1$ at $(j,i)$ and 0's at all the other positions, this would be 
the isomorphism of $\sspan\{f_{ij}\}$ and $\gotso_n$.\\
Hence in this case $ \I_2 \simeq \gotsl_2 \oplus \gotso_n$.

\subsection{Twisted cubic}
Now starting from the exceptional group $G_2$ one gets to the example of the 
twisted cubic. 

\begin{figure}[htb]
\center
\includegraphics[width=0.9\textwidth]{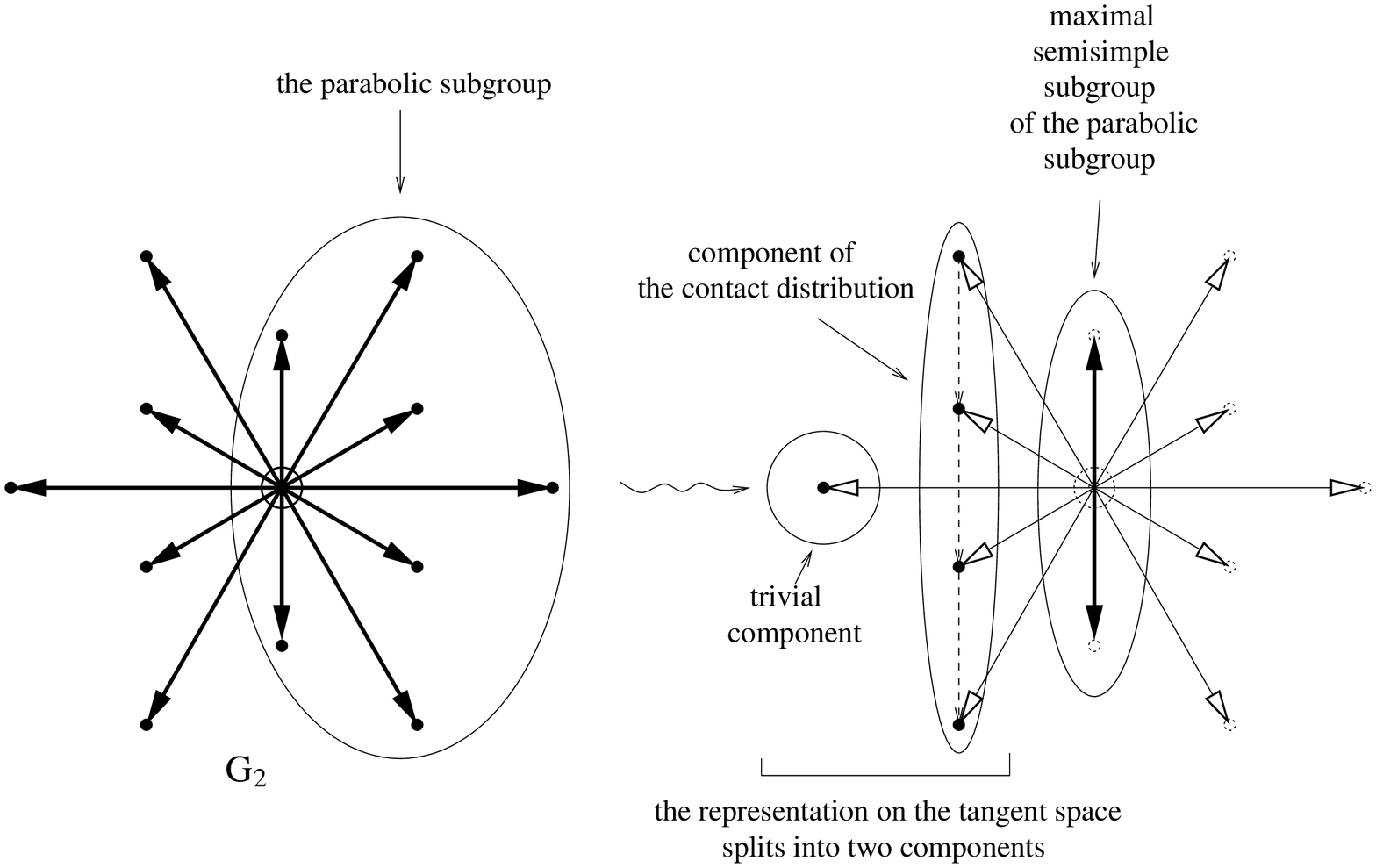}
\caption{The weight system of the adjoint action of $G_2$. When restricted to the maximal semisimple
subgroup of the parabolic subgroup it decomposes into several components. In particular the tangent space to
the minimal nilpotent orbit splits into two components: a trivial one and the fibre of the contact distribution.
The last one in this case is the third symmetric power of the standard representation of $\gotsl_2$  and the closed
orbit is the twisted cubic.}
 \label{pic_adjoint_G2}
\end{figure}

Let me have one more look on the example \ref{przyklad_skr_kubika} and examine it's ideal.\\

 $\I$ is generated by functions $f_+, f_-, h' \in \C[x_0,x_1,x_2,x_3]$:\nopagebreak
\begin{itemize}
  \item [] $f_+:= x_2^2 - x_1x_3$;
  \item [] $f_-:= x_0x_2 - x_1^2$;
  \item [] $h':= x_0x_3 - x_1x_2$.
\end{itemize}
Their derivations are:
\begin{itemize}
  \item [] $\ud f_+= (  0,  -x_3, 2x_2, -x_1)$;
  \item [] $\ud f_-= (x_2, -2x_1,  x_0,    0)$;
  \item [] $\ud h'=  (x_3,  -x_2, -x_1,  x_0)$.
\end{itemize}
and the Lie brackets (recall that the matrix of $\omega'$ is now
$\begin{pmatrix}
  0 & 0 & 0 & 3 \\
  0 & 0 &-1 & 0 \\
  0 & 1 & 0 & 0 \\
 -3 & 0 & 0 & 0
  \end{pmatrix}$):\nopagebreak
\begin{itemize}
  \item [] $[f_+,f_-] = x_3x_0 -4x_2x_1 +3x_1 x_2 = x_0x_3 - x_1x_2 = h'$;
  \item [] $[h' ,f_+] = -3x_3x_1 + 2x_2^2 + x_1x_3 = 2 (x_2^2 - x_1x_3) = 2 f_+$;
  \item [] $[h' ,f_-] = x_2x_0 + 2x_1^2 -3x_0x_2 = -2(x_2x_0 - x_1^2) = -2 f_-$;
\end{itemize}

So $\I_2$ in this case is isomorphic to $\gotsl_2$.

\subsection{Grassmannian $Gr(3,6)$} \label{example_grassmannian}
The adjoint variety for the exceptional group $E_6$ gives rise to legendrian variety $Gr(3,6)$  - the full 
Grassmannian of $3-$spaces in a $6-$space.

\begin{figure}[htb]
\centering
\includegraphics[width=0.9\textwidth]{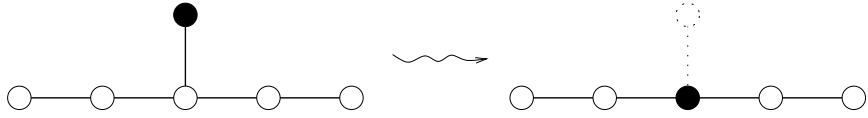}
\caption{The two homogeneous varieties for the group $E_6$. }
 \label{pic_adjoint_E6}
\end{figure}

The Grassmannian $Gr(3,6)$ is naturally embedded in the projectivisation of $V:=\Lambda^3 \C^6$:
 the embedding is given by assigning $[u \wedge v \wedge w]$ to the space spanned by $\{u,v,w\}$.
The symplectic form on $V$ is given by the natural map:
$$\omega: \ \Lambda^2 V = \Lambda^2(\Lambda^3 \C^6)  \lra \Lambda^6 \C^6 \simeq \C^6.$$

Further let $a_1 \ldots a_6$ be some coordinates on $\C^6$ and let $x_{ijk}:=a_i \wedge a_j \wedge a_k$. Then the 35 equations of
$Gr(3,6)$ are:
$$ \begin{array}{cc}
    x_{243} x_{456} + x_{346} x_{245} - x_{354} x_{264} = 0,&
    x_{235} x_{456} + x_{365} x_{245} - x_{354} x_{256} = 0,\\
    x_{263} x_{245} - x_{235} x_{264} + x_{243} x_{256} = 0,&
    x_{263} x_{354} - x_{235} x_{346} + x_{243} x_{365} = 0,\\
    x_{263} x_{456} + x_{365} x_{264} - x_{346} x_{256} = 0,&
    x_{125} x_{354} + x_{135} x_{245} - x_{235} x_{145} = 0,\\
    x_{124} x_{256} + x_{125} x_{264} + x_{126} x_{245} = 0,&
    x_{124} x_{354} + x_{134} x_{245} + x_{243} x_{145} = 0,\\
    x_{124} x_{346} + x_{134} x_{264} - x_{243} x_{146} = 0,&
    -x_{123} x_{245} + x_{124} x_{235} + x_{125} x_{243} = 0,\\
    -x_{123} x_{264} + x_{124} x_{263} - x_{126} x_{243} = 0,&
    x_{123} x_{256} + x_{125} x_{263} + x_{126} x_{235} = 0,\\
    x_{123} x_{354} + x_{134} x_{235} + x_{135} x_{243} = 0,&
    x_{123} x_{346} + x_{134} x_{263} - x_{136} x_{243} = 0,\\
    \multicolumn{2}{c}{
    -x_{123} x_{365} + x_{135} x_{263} + x_{136} x_{235} = 0,}
\end{array}$$
$$ \begin{array}{cc}
    x_{136} x_{456} + x_{365} x_{146} + x_{346} x_{156} = 0,&
    x_{135} x_{456} + x_{365} x_{145} - x_{354} x_{156} = 0,\\
    x_{134} x_{156} - x_{135} x_{146} + x_{136} x_{145} = 0,&
    x_{134} x_{365} + x_{135} x_{346} + x_{136} x_{354} = 0,\\
    x_{134} x_{456} - x_{346} x_{145} - x_{354} x_{146} = 0,&
    x_{126} x_{346} + x_{136} x_{264} - x_{263} x_{146} = 0,\\
    x_{126} x_{365} + x_{136} x_{256} + x_{263} x_{156} = 0,&
    x_{126} x_{456} - x_{256} x_{146} - x_{264} x_{156} = 0,\\
    x_{125} x_{365} + x_{135} x_{256} - x_{235} x_{156} = 0,&
    x_{125} x_{456} - x_{256} x_{145} + x_{245} x_{156} = 0,\\
    x_{124} x_{156} - x_{125} x_{146} + x_{126} x_{145} = 0,&
    x_{124} x_{456} + x_{264} x_{145} + x_{245} x_{146} = 0,\\
    x_{123} x_{156} - x_{125} x_{136} + x_{126} x_{135} = 0,&
    x_{123} x_{146} - x_{124} x_{136} + x_{126} x_{134} = 0,\\
    \multicolumn{2}{c}{
    x_{123} x_{145} - x_{124} x_{135} + x_{125} x_{134} = 0,}
\end{array}$$
$$ \begin{array}{c}
    x_{126} x_{354} + x_{136} x_{245} - x_{235} x_{146} - x_{243} x_{156} =0,\\
    x_{125} x_{346} + x_{135} x_{264} - x_{263} x_{145} - x_{243} x_{156} =0,\\
    x_{124} x_{365} - x_{135} x_{264} - x_{136} x_{245} + x_{243} x_{156} =0,\\
    x_{123} x_{456} + x_{263} x_{145} + x_{235} x_{146} + x_{243} x_{156} =0,\\
x_{134} x_{256} + x_{135} x_{264} + x_{136} x_{245} - x_{263} x_{145} - x_{235} x_{146} - x_{243} x_{156} = 0;
\end{array} $$

The Lie algebra $\gotg$  in this case turns out to be $\gotsl_6$ (not really surprising \ldots).
The first 15 polynomials correspond to the positive roots of $\gotsl_6$, next 15 to the negative roots while the last 5 of them span the Cartan subalgebra.

\subsection{Lagrangian Grassmannian $Gr_L(3,6)$} \label{example_lagrangian_grassmannian}
The next example arises from the adjoint variety for the exceptional group $F_4$:
it is a legendrian variety $Gr_L(3,6)$  - the Grassmannian of Lagrangian  $3-$spaces in a $6-$space.

\begin{figure}[htb]
\centering
\includegraphics[width=0.9\textwidth]{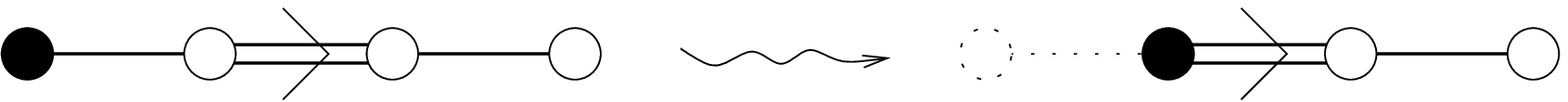}
\caption{The two homogeneous varieties for the group $F_4$. }
 \label{pic_adjoint_F4}
\end{figure}

Assume on $\C^6$ a symplectic form $\omega_{\C^6}$ is given. Then it determines a map:
$$
\begin{array}{rcl}
\Lambda^3 \C^6      & \lra        & \Lambda^1 \C^6 \simeq \C^6 \\
u \wedge v \wedge w & \longmapsto & \omega_{\C^6}(u,v) w + \omega_{\C^6}(v,w) u + \omega_{\C^6}(w,u) v
\end{array}
$$
kernel of which is a 14 dimensional subspace $V \subset \Lambda^3 \C^6$. 
Intersecting $Gr(3,6) \cap \P(V)$ one gets the Grassmannian of Lagrangian subspaces in $\C^6$. The equations of $V$ are:
$$
\begin{array}{ccc}
x_{124} + x_{263}=0,&
x_{125} + x_{136}=0,&
x_{134} + x_{235}=0,\\
x_{365} + x_{145}=0,&
x_{346} + x_{245}=0,&
x_{256} + x_{146}=0;
\end{array}
$$
and the 35 equations of $Gr(3,6)$ reduce to 21 equations of $Gr_L(3,6)$: 
$$ \begin{array}{cc}
    4 y_{3} y_{7} - 4 y_{8} y_{9} - y_{12}^2 ,&
    -4 y_{6} y_{7} - 2 y_{8} y_{13} - y_{11} y_{12} ,\\
    y_{3} y_{13} + y_{4} y_{12} + 2 y_{6} y_{9} ,&
    y_{3} y_{11} - 2 y_{4} y_{8} + y_{6} y_{12} ,\\
    -4 y_{4} y_{7} + 2 y_{9} y_{11} - y_{12} y_{13} ,&
    -y_{1} y_{12} + y_{4} y_{13} + 2 y_{5} y_{9} ,\\
    y_{0} y_{12} + 2 y_{3} y_{5} - 2 y_{4} y_{6} ,&
    -y_{0} y_{9} - y_{1} y_{3} - y_{4}^2 ,\\
\multicolumn{2}{c}{   
    y_{0} y_{8} + y_{2} y_{3} - y_{6}^2,}
\end{array}$$
$$ \begin{array}{cc}
    -y_{2} y_{12} + 2 y_{5} y_{8} - y_{6} y_{11} ,&
    y_{0} y_{13} - 2 y_{1} y_{6} - 2 y_{4} y_{5} ,\\
    -y_{0} y_{11} - 2 y_{2} y_{4} + 2 y_{5} y_{6} ,&
    -4 y_{5} y_{7} + 2 y_{10} y_{12} - y_{11} y_{13} ,\\
    4 y_{2} y_{7} - 4 y_{8} y_{10} - y_{11}^2 ,&
    y_{2} y_{13} + y_{5} y_{11} + 2 y_{6} y_{10} ,\\
    y_{1} y_{11} - 2 y_{4} y_{10} - y_{5} y_{13} ,&
    4 y_{1} y_{7} - 4 y_{9} y_{10} + y_{13}^2 ,\\
\multicolumn{2}{c}{
    y_{0} y_{10} + y_{1} y_{2} + y_{5}^2 ,}
\end{array}$$
$$ \begin{array}{c}
    2 y_{2} y_{9} - 2 y_{3} y_{10} - y_{4} y_{11} + y_{5} y_{12} ,\\
    2 y_{1} y_{8} - 2 y_{3} y_{10} + y_{5} y_{12} - y_{6} y_{13} ,\\
    2 y_{0} y_{7} + 2 y_{3} y_{10} + y_{4} y_{11} + y_{6} y_{13};
\end{array}$$
where $y_i$'s are the following symplectic coordinates:
$$
\begin{array}{cc}
y_{0}:= x_{123}, &
y_{7}:= x_{456},\\
y_{1}:= x_{126}, &
y_{8}:= x_{354},\\
y_{2}:= x_{135}, &
y_{9}:= x_{264},\\
y_{3}:= x_{243}, &
y_{10}:= x_{156},\\
y_{4}:= x_{124} = -x_{263},&
y_{11}:= 2x_{365} = -2x_{145},\\
y_{5}:= x_{125} = -x_{136},&
y_{12}:= 2x_{346}= -2x_{245},\\
y_{6}:= x_{134} = -x_{235}, &
y_{13}:= 2x_{256}= -2x_{146};
\end{array}
$$

The Lie algebra $\gotg$ is now isomorphic to $\gotsp_6$.

\subsection{Spinor variety $\mathbb{S}_6$} \label{example_spinor}
Out of the adjoint variety for the exceptional group $E_7$ one gets as the 
legendrian variety the spinor variety $\mathbb{S}_6$.

\begin{figure}[htb]
\centering
\includegraphics[width=0.9\textwidth]{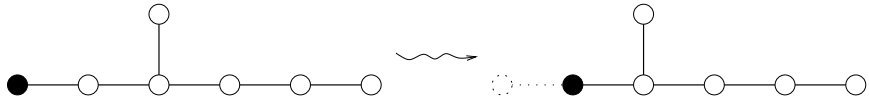}
\caption{The two homogeneous varieties for the group $E_7$. }
 \label{pic_adjoint_E7}
\end{figure}

The spinor variety  $\mathbb{S}_6$ parametrises all the linear $\P^5$'s on a smooth quadric $Q^{10}\subset \P^{11}$. 
It can be embedded in $\P(V)$ for 
$$V:= \Lambda^{\textrm{even}}\C^6 = \Lambda^0 \C^6 \osum \Lambda^2 \C^6 \osum \Lambda^4 \C^6 \osum \Lambda^6 \C^6.$$
Define the coordinates on $V$ as follows:
$$
\begin{array}{rl}
x                         & \textrm{the coordinate  on }  \Lambda^0 \C^6 \\
m_{ij} = a_i \wedge a_j   & \textrm{the coordinates on }  \Lambda^2 \C^6 \\
n_{ij} = 
a_1 \wedge \stackrel {\stackrel {i,j} \lor} \ldots \wedge a_6  
 & \textrm{the coordinates on }  \Lambda^4 \C^6 \\
y = a_1 \wedge \ldots \wedge a_6  &\textrm{the coordinate on }  \Lambda^6 \C^6 
\end{array}
$$
Now let $M$ and $N$ be the skew-symmetric matrices of $m_{ij}$'s and $n_{ij}$'s correspondingly. 
Also let $Pf_M$ and $Pf_N$ be the matrices of all $4\times 4$ Pfaffians of $M$ and $N$.
Then the 66 equations of the spinor variety are of the following three types:
$$
M N = x y \Id_6, \quad Pf_M = - x N, \quad Pf_N = y M
$$ 
(more detailed treatment of a spinor variety (but $\mathbb{S}_5$) is in \cite{reid_corti}).
The Lie algebra $\gotg$ for $\mathbb{S}_6$ is $\gotso_{12}$.

\subsection{$E_7$ variety } \label{example_E7}
Finally, the last example arises from the adjoint variety of $E_8$.
\begin{figure}[htb]
\centering
\includegraphics[width=0.5\textwidth]{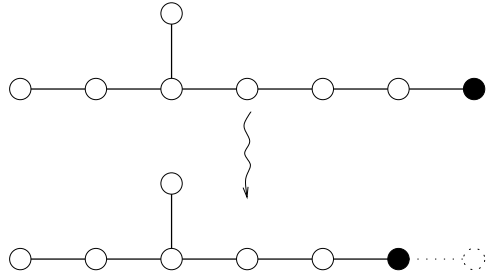}
\caption{The two homogeneous varieties for the group $E_8$. }
 \label{pic_adjoint_E8}
\end{figure}
It is a 27-dimensional $E_7$ variety embedded in $\P^{55}$. The 133 equations are:
$$
 \begin{array}{cc}
\scriptstyle{x_{7} x_{34} + x_{9} x_{36} + x_{11} x_{38} + x_{14} x_{40} + x_{17} x_{43} + x_{27} x_{54}=0,}&
\scriptstyle{x_{5} x_{32} + x_{8} x_{34} + x_{9} x_{35} + x_{20} x_{46} +  x_{22} x_{49} + x_{24} x_{51}=0,} \\ 
\scriptstyle{
    x_{6} x_{32} + x_{8} x_{33} + x_{13} x_{39} + x_{16} x_{42} +  x_{19} x_{45} + x_{26} x_{53}=0, }& \scriptstyle{
    x_{4} x_{31} + x_{10} x_{36} + x_{11} x_{37} + x_{18} x_{44} +  x_{21} x_{47} + x_{25} x_{52}=0, }\\  \scriptstyle{
    x_{3} x_{30} + x_{12} x_{38} + x_{14} x_{39} + x_{16} x_{41} +  x_{23} x_{49} + x_{24} x_{50}=0, }&  \scriptstyle{
    x_{2} x_{29} + x_{15} x_{40} + x_{17} x_{42} + x_{19} x_{44} +  x_{21} x_{46} + x_{22} x_{48}=0, }\\  \scriptstyle{
    x_{1} x_{28} + x_{43} x_{55} - x_{45} x_{54} + x_{47} x_{53} - x_{49} x_{52} + x_{50} x_{51}=0, }& 
    \scriptstyle{
    x_{7} x_{32} + x_{9} x_{33} - x_{13} x_{38} - x_{16} x_{40} -  x_{19} x_{43} + x_{27} x_{53}=0, }\\  \scriptstyle{
    x_{5} x_{31} - x_{10} x_{34} - x_{11} x_{35} + x_{20} x_{44} +  x_{22} x_{47} - x_{25} x_{51}=0, }&  \scriptstyle{
    x_{6} x_{31} - x_{10} x_{33} + x_{13} x_{37} - x_{18} x_{42} -   x_{21} x_{45} + x_{26} x_{52}=0, }\\  \scriptstyle{
    x_{4} x_{30} - x_{12} x_{36} - x_{14} x_{37} + x_{18} x_{41} -   x_{23} x_{47} + x_{25} x_{50}=0, }&  \scriptstyle{
    x_{3} x_{29} - x_{15} x_{38} - x_{17} x_{39} - x_{19} x_{41} +   x_{23} x_{46} + x_{24} x_{48}=0, }\\  \scriptstyle{
    x_{2} x_{28} - x_{40} x_{55} + x_{42} x_{54} - x_{44} x_{53} +   x_{46} x_{52} - x_{48} x_{51}=0, }&  \scriptstyle{
    x_{7} x_{31} - x_{11} x_{33} - x_{13} x_{36} + x_{18} x_{40} +   x_{21} x_{43} + x_{27} x_{52}=0, }\\  \scriptstyle{
   x_{8} x_{31} + x_{10} x_{32} - x_{13} x_{35} - x_{20} x_{42} -   x_{22} x_{45} - x_{26} x_{51}=0, }&
\scriptstyle{
    x_{5} x_{30} + x_{12} x_{34} + x_{14} x_{35} + x_{20} x_{41} -   x_{24} x_{47} - x_{25} x_{49}=0, }
\end{array}
$$
$$
 \begin{array}{cc}
  \scriptstyle{
    x_{6} x_{30} + x_{12} x_{33} - x_{16} x_{37} - x_{18} x_{39} +   x_{23} x_{45} + x_{26} x_{50}=0, }&  \scriptstyle{
    x_{4} x_{29} + x_{15} x_{36} + x_{17} x_{37} - x_{21} x_{41} -   x_{23} x_{44} + x_{25} x_{48}=0, }\\  \scriptstyle{
    x_{3} x_{28} + x_{38} x_{55} - x_{39} x_{54} + x_{41} x_{53} -   x_{46} x_{50} + x_{48} x_{49}=0, }&  \scriptstyle{
    x_{9} x_{31} + x_{11} x_{32} + x_{13} x_{34} + x_{20} x_{40} +   x_{22} x_{43} - x_{27} x_{51}=0, }\\  \scriptstyle{
    x_{7} x_{30} + x_{14} x_{33} + x_{16} x_{36} + x_{18} x_{38} -   x_{23} x_{43} + x_{27} x_{50}=0, }&  \scriptstyle{
    x_{8} x_{30} - x_{12} x_{32} + x_{16} x_{35} - x_{20} x_{39} +   x_{24} x_{45} - x_{26} x_{49}=0, }\\  \scriptstyle{
    x_{5} x_{29} - x_{15} x_{34} - x_{17} x_{35} - x_{22} x_{41} -   x_{24} x_{44} - x_{25} x_{46}=0, }&\scriptstyle{
    x_{6} x_{29} - x_{15} x_{33} + x_{19} x_{37} + x_{21} x_{39} +   x_{23} x_{42} + x_{26} x_{48}=0, }\\  \scriptstyle{
    x_{4} x_{28} - x_{36} x_{55} + x_{37} x_{54} - x_{41} x_{52} +   x_{44} x_{50} - x_{47} x_{48}=0, }&  \scriptstyle{
    x_{9} x_{30} - x_{14} x_{32} - x_{16} x_{34} + x_{20} x_{38} -   x_{24} x_{43} - x_{27} x_{49}=0, }\\  \scriptstyle{
    x_{7} x_{29} - x_{17} x_{33} - x_{19} x_{36} - x_{21} x_{38} -   x_{23} x_{40} + x_{27} x_{48}=0, }&  \scriptstyle{
    x_{10} x_{30} + x_{12} x_{31} + x_{18} x_{35} + x_{20} x_{37} +  x_{25} x_{45} + x_{26} x_{47}=0, }\\  \scriptstyle{
    x_{8} x_{29} + x_{15} x_{32} - x_{19} x_{35} + x_{22} x_{39} +   x_{24} x_{42} - x_{26} x_{46}=0, }&  \scriptstyle{
    x_{5} x_{28} + x_{34} x_{55} - x_{35} x_{54} + x_{41} x_{51} -   x_{44} x_{49} + x_{46} x_{47}=0, }\\  \scriptstyle{
    x_{6} x_{28} + x_{33} x_{55} - x_{37} x_{53} + x_{39} x_{52} -   x_{42} x_{50} + x_{45} x_{48}=0, }&
  \scriptstyle{
    x_{11} x_{30} + x_{14} x_{31} - x_{18} x_{34} - x_{20} x_{36} -  x_{25} x_{43} + x_{27} x_{47}=0, }
\end{array}
$$
$$
 \begin{array}{cc}
  \scriptstyle{
    x_{9} x_{29} + x_{17} x_{32} + x_{19} x_{34} - x_{22} x_{38} -   x_{24} x_{40} - x_{27} x_{46}=0, }&  \scriptstyle{
    x_{7} x_{28} - x_{33} x_{54} + x_{36} x_{53} - x_{38} x_{52} +   x_{40} x_{50} - x_{43} x_{48}=0, }\\  \scriptstyle{
    x_{10} x_{29} - x_{15} x_{31} - x_{21} x_{35} - x_{22} x_{37} +  x_{25} x_{42} + x_{26} x_{44}=0, }&  \scriptstyle{
    x_{8} x_{28} - x_{32} x_{55} + x_{35} x_{53} - x_{39} x_{51} +   x_{42} x_{49} - x_{45} x_{46}=0, }\\  \scriptstyle{
    x_{13} x_{30} + x_{16} x_{31} + x_{18} x_{32} + x_{20} x_{33} -  x_{26} x_{43} - x_{27} x_{45}=0, }&  \scriptstyle{
    x_{11} x_{29} - x_{17} x_{31} + x_{21} x_{34} + x_{22} x_{36} -  x_{25} x_{40} + x_{27} x_{44}=0, }\\  \scriptstyle{
    x_{9} x_{28} + x_{32} x_{54} - x_{34} x_{53} + x_{38} x_{51} -   x_{40} x_{49} + x_{43} x_{46}=0, }&
  \scriptstyle{
    x_{12} x_{29} + x_{15} x_{30} - x_{23} x_{35} - x_{24} x_{37} -  x_{25} x_{39} - x_{26} x_{41}=0, }\\  \scriptstyle{
    x_{10} x_{28} + x_{31} x_{55} - x_{35} x_{52} + x_{37} x_{51} -  x_{42} x_{47} + x_{44} x_{45}=0, }&   \scriptstyle{
    x_{13} x_{29} - x_{19} x_{31} - x_{21} x_{32} - x_{22} x_{33} -  x_{26} x_{40} - x_{27} x_{42}=0, }\\  \scriptstyle{
    x_{14} x_{29} + x_{17} x_{30} + x_{23} x_{34} + x_{24} x_{36} +  x_{25} x_{38} - x_{27} x_{41}=0, }&   \scriptstyle{
    x_{11} x_{28} - x_{31} x_{54} + x_{34} x_{52} - x_{36} x_{51} +  x_{40} x_{47} - x_{43} x_{44}=0, }\\  \scriptstyle{
    x_{12} x_{28} - x_{30} x_{55} + x_{35} x_{50} - x_{37} x_{49} +  x_{39} x_{47} - x_{41} x_{45}=0, } &  \scriptstyle{
    x_{16} x_{29} + x_{19} x_{30} - x_{23} x_{32} - x_{24} x_{33} +  x_{26} x_{38} + x_{27} x_{39}=0, }\\  \scriptstyle{
    x_{13} x_{28} + x_{31} x_{53} - x_{32} x_{52} + x_{33} x_{51} -  x_{40} x_{45} + x_{42} x_{43}=0, }&   \scriptstyle{
    x_{14} x_{28} + x_{30} x_{54} - x_{34} x_{50} + x_{36} x_{49} -  x_{38} x_{47} + x_{41} x_{43}=0, }
\end{array}
$$
$$
 \begin{array}{cc}
\scriptstyle{
    x_{15} {28} + x_{29} x_{55} - x_{35} x_{48} + x_{37} x_{46} -  x_{39} x_{44} + x_{41} x_{42}=0, }&   \scriptstyle{
    x_{18} x_{29} + x_{21} x_{30} + x_{23} x_{31} - x_{25} x_{33} -  x_{26} x_{36} - x_{27} x_{37}=0, }\\  \scriptstyle{
    x_{16} x_{28} - x_{30} x_{53} + x_{32} x_{50} - x_{33} x_{49} +  x_{38} x_{45} - x_{39} x_{43}=0, }&   \scriptstyle{
    x_{17} x_{28} - x_{29} x_{54} + x_{34} x_{48} - x_{36} x_{46} +  x_{38} x_{44} - x_{40} x_{41}=0, }\\  \scriptstyle{
    x_{20} x_{29} + x_{22} x_{30} + x_{24} x_{31} + x_{25} x_{32} +  x_{26} x_{34} + x_{27} x_{35}=0, } &  \scriptstyle{
    x_{18} x_{28} + x_{30} x_{52} - x_{31} x_{50} + x_{33} x_{47} -  x_{36} x_{45} + x_{37} x_{43}=0, }\\  \scriptstyle{
    x_{19} x_{28} + x_{29} x_{53} - x_{32} x_{48} + x_{33} x_{46} -  x_{38} x_{42} + x_{39} x_{40}=0, } &  \scriptstyle{
    x_{20} x_{28} - x_{30} x_{51} + x_{31} x_{49} - x_{32} x_{47} +  x_{34} x_{45} - x_{35} x_{43}=0, }\\  \scriptstyle{
    x_{21} x_{28} - x_{29} x_{52} + x_{31} x_{48} - x_{33} x_{44} +  x_{36} x_{42} - x_{37} x_{40}=0, }&   \scriptstyle{
    x_{22} x_{28} + x_{29} x_{51} - x_{31} x_{46} + x_{32} x_{44} -  x_{34} x_{42} + x_{35} x_{40}=0, }\\  \scriptstyle{
    x_{23} x_{28} + x_{29} x_{50} - x_{30} x_{48} + x_{33} x_{41} -  x_{36} x_{39} + x_{37} x_{38}=0, }&   \scriptstyle{
    x_{24} x_{28} - x_{29} x_{49} + x_{30} x_{46} - x_{32} x_{41} +  x_{34} x_{39} - x_{35} x_{38}=0, }\\  \scriptstyle{
    x_{25} x_{28} + x_{29} x_{47} - x_{30} x_{44} + x_{31} x_{41} -  x_{34} x_{37} + x_{35} x_{36}=0, }&   \scriptstyle{
    x_{26} x_{28} - x_{29} x_{45} + x_{30} x_{42} - x_{31} x_{39} +  x_{32} x_{37} - x_{33} x_{35}=0, }\\  \scriptstyle{
    x_{27} x_{28} + x_{29} x_{43} - x_{30} x_{40} + x_{31} x_{38} -  x_{32} x_{36} + x_{33} x_{34}=0, }
\end{array}
$$
$$
 \begin{array}{cc}
\scriptstyle{
    x_{6} x_{35} + x_{8} x_{37} + x_{10} x_{39} + x_{12} x_{42} +   x_{15} x_{45} + x_{26} x_{55}=0, }&  \scriptstyle{
    x_{4} x_{33} + x_{6} x_{36} + x_{7} x_{37} + x_{18} x_{48} + x_{21} x_{50} + x_{23} x_{52}=0, }\\   \scriptstyle{
    x_{4} x_{34} + x_{5} x_{36} + x_{11} x_{41} + x_{14} x_{44} + x_{17} x_{47} + x_{25} x_{54}=0, }&  \scriptstyle{
    x_{3} x_{32} + x_{8} x_{38} + x_{9} x_{39} + x_{16} x_{46} +  x_{19} x_{49} + x_{24} x_{53}=0, }\\   \scriptstyle{
    x_{2} x_{31} + x_{10} x_{40} + x_{11} x_{42} + x_{13} x_{44} + x_{21} x_{51} + x_{22} x_{52}=0, }&  \scriptstyle{
    x_{1} x_{30} + x_{12} x_{43} + x_{14} x_{45} + x_{16} x_{47} + x_{18} x_{49} + x_{20} x_{50}=0, }\\   \scriptstyle{
    x_{0} x_{29} - x_{15} x_{27} + x_{17} x_{26} - x_{19} x_{25} + x_{21} x_{24} - x_{22} x_{23}=0, }&  \scriptstyle{
    x_{4} x_{35} + x_{5} x_{37} - x_{10} x_{41} - x_{12} x_{44} -  x_{15} x_{47} + x_{25} x_{55}=0, }\\   \scriptstyle{
    x_{3} x_{33} - x_{6} x_{38} - x_{7} x_{39} + x_{16} x_{48} +   x_{19} x_{50} - x_{23} x_{53}=0, }&  \scriptstyle{
    x_{3} x_{34} - x_{5} x_{38} + x_{9} x_{41} - x_{14} x_{46} -   x_{17} x_{49} + x_{24} x_{54}=0, }\\   \scriptstyle{
    x_{2} x_{32} - x_{8} x_{40} - x_{9} x_{42} + x_{13} x_{46} -   x_{19} x_{51} + x_{22} x_{53}=0, }&  \scriptstyle{
    x_{1} x_{31} - x_{10} x_{43} - x_{11} x_{45} - x_{13} x_{47} +  x_{18} x_{51} + x_{20} x_{52}=0, }\\   \scriptstyle{
    x_{0} x_{30} + x_{12} x_{27} - x_{14} x_{26} + x_{16} x_{25} -  x_{18} x_{24} + x_{20} x_{23}=0, }&  \scriptstyle{
    x_{3} x_{35} - x_{5} x_{39} - x_{8} x_{41} + x_{12} x_{46} +    x_{15} x_{49} + x_{24} x_{55}=0, }\\   \scriptstyle{
    x_{3} x_{36} + x_{4} x_{38} - x_{7} x_{41} - x_{14} x_{48} -    x_{17} x_{50} - x_{23} x_{54}=0, }&  \scriptstyle{
    x_{2} x_{33} + x_{6} x_{40} + x_{7} x_{42} + x_{13} x_{48} -    x_{19} x_{52} - x_{21} x_{53}=0, }
\end{array}
$$
$$
 \begin{array}{cc}
   \scriptstyle{
    x_{2} x_{34} + x_{5} x_{40} - x_{9} x_{44} - x_{11} x_{46} +    x_{17} x_{51} + x_{22} x_{54}=0, }&  \scriptstyle{
    x_{1} x_{32} + x_{8} x_{43} + x_{9} x_{45} - x_{13} x_{49} -    x_{16} x_{51} + x_{20} x_{53}=0, }\\   \scriptstyle{
    x_{0} x_{31} - x_{10} x_{27} + x_{11} x_{26} - x_{13} x_{25} +  x_{18} x_{22} - x_{20} x_{21}=0, }&  \scriptstyle{
    x_{3} x_{37} + x_{4} x_{39} + x_{6} x_{41} + x_{12} x_{48} +    x_{15} x_{50} - x_{23} x_{55}=0, }\\   \scriptstyle{
    x_{2} x_{35} + x_{5} x_{42} + x_{8} x_{44} + x_{10} x_{46} -    x_{15} x_{51} + x_{22} x_{55}=0, }&  \scriptstyle{
    x_{2} x_{36} - x_{4} x_{40} + x_{7} x_{44} - x_{11} x_{48} +    x_{17} x_{52} - x_{21} x_{54}=0, }\\   \scriptstyle{
    x_{1} x_{33} - x_{6} x_{43} - x_{7} x_{45} - x_{13} x_{50} -    x_{16} x_{52} - x_{18} x_{53}=0, }&  \scriptstyle{
    x_{1} x_{34} - x_{5} x_{43} + x_{9} x_{47} + x_{11} x_{49} +    x_{14} x_{51} + x_{20} x_{54}=0, }\\   \scriptstyle{
    x_{0} x_{32} + x_{8} x_{27} - x_{9} x_{26} + x_{13} x_{24} -   x_{16} x_{22} + x_{19} x_{20}=0, }&  \scriptstyle{
    x_{2} x_{37} - x_{4} x_{42} - x_{6} x_{44} + x_{10} x_{48} -   x_{15} x_{52} - x_{21} x_{55}=0, }\\   \scriptstyle{
    x_{1} x_{35} - x_{5} x_{45} - x_{8} x_{47} - x_{10} x_{49} -   x_{12} x_{51} + x_{20} x_{55}=0, }&  \scriptstyle{
    x_{2} x_{38} + x_{3} x_{40} + x_{7} x_{46} + x_{9} x_{48} +    x_{17} x_{53} + x_{19} x_{54}=0, }\\   \scriptstyle{
    x_{1} x_{36} + x_{4} x_{43} - x_{7} x_{47} + x_{11} x_{50} +    x_{14} x_{52} - x_{18} x_{54}=0, }&  \scriptstyle{
    x_{0} x_{33} - x_{6} x_{27} + x_{7} x_{26} - x_{13} x_{23} +    x_{16} x_{21} - x_{18} x_{19}=0, }\\   \scriptstyle{
    x_{0} x_{34} - x_{5} x_{27} + x_{9} x_{25} - x_{11} x_{24} +    x_{14} x_{22} - x_{17} x_{20}=0, }&  \scriptstyle{
    x_{2} x_{39} + x_{3} x_{42} - x_{6} x_{46} - x_{8} x_{48} -     x_{15} x_{53} + x_{19} x_{55}=0, }
\end{array}
$$
$$
 \begin{array}{cc}
   \scriptstyle{
    x_{1} x_{37} + x_{4} x_{45} + x_{6} x_{47} - x_{10} x_{50} -    x_{12} x_{52} - x_{18} x_{55}=0, }&  \scriptstyle{
    x_{0} x_{35} + x_{5} x_{26} - x_{8} x_{25} + x_{10} x_{24} -    x_{12} x_{22} + x_{15} x_{20}=0, }\\   \scriptstyle{
    x_{1} x_{38} - x_{3} x_{43} - x_{7} x_{49} - x_{9} x_{50} +     x_{14} x_{53} + x_{16} x_{54}=0, }&  \scriptstyle{
    x_{0} x_{36} + x_{4} x_{27} - x_{7} x_{25} + x_{11} x_{23} -    x_{14} x_{21} + x_{17} x_{18}=0, }\\   \scriptstyle{
    x_{2} x_{41} + x_{3} x_{44} + x_{4} x_{46} + x_{5} x_{48} -     x_{15} x_{54} - x_{17} x_{55}=0, }&  \scriptstyle{
    x_{1} x_{39} - x_{3} x_{45} + x_{6} x_{49} + x_{8} x_{50} -     x_{12} x_{53} + x_{16} x_{55}=0, }\\   \scriptstyle{
    x_{0} x_{37} - x_{4} x_{26} + x_{6} x_{25} - x_{10} x_{23} +    x_{12} x_{21} - x_{15} x_{18}=0, }&  \scriptstyle{
    x_{1} x_{40} + x_{2} x_{43} - x_{7} x_{51} - x_{9} x_{52} -     x_{11} x_{53} - x_{13} x_{54}=0, }\\   \scriptstyle{
    x_{0} x_{38} - x_{3} x_{27} + x_{7} x_{24} - x_{9} x_{23} +     x_{14} x_{19} - x_{16} x_{17}=0, }&  \scriptstyle{
    x_{1} x_{41} - x_{3} x_{47} - x_{4} x_{49} - x_{5} x_{50} -     x_{12} x_{54} - x_{14} x_{55}=0, }\\   \scriptstyle{
    x_{1} x_{42} + x_{2} x_{45} + x_{6} x_{51} + x_{8} x_{52} +     x_{10} x_{53} - x_{13} x_{55}=0, }&  \scriptstyle{
    x_{0} x_{39} + x_{3} x_{26} - x_{6} x_{24} + x_{8} x_{23} -     x_{12} x_{19} + x_{15} x_{16}=0, }\\   \scriptstyle{
    x_{0} x_{40} + x_{2} x_{27} - x_{7} x_{22} + x_{9} x_{21} -     x_{11} x_{19} + x_{13} x_{17}=0, }&  \scriptstyle{
    x_{1} x_{44} + x_{2} x_{47} - x_{4} x_{51} - x_{5} x_{52} +     x_{10} x_{54} + x_{11} x_{55}=0, }\\   \scriptstyle{
    x_{0} x_{41} - x_{3} x_{25} + x_{4} x_{24} - x_{5} x_{23} +     x_{12} x_{17} - x_{14} x_{15}=0, }&  \scriptstyle{
    x_{0} x_{42} - x_{2} x_{26} + x_{6} x_{22} - x_{8} x_{21} +     x_{10} x_{19} - x_{13} x_{15}=0, }
\end{array}
$$
$$
 \begin{array}{cc}
   \scriptstyle{
    x_{0} x_{43} - x_{1} x_{27} + x_{7} x_{20} - x_{9} x_{18} +    x_{11} x_{16} - x_{13} x_{14}=0, }&  \scriptstyle{
    x_{1} x_{46} + x_{2} x_{49} + x_{3} x_{51} - x_{5} x_{53} -     x_{8} x_{54} - x_{9} x_{55}=0, }\\   \scriptstyle{
    x_{0} x_{44} + x_{2} x_{25} - x_{4} x_{22} + x_{5} x_{21} -      x_{10} x_{17} + x_{11} x_{15}=0, }&  \scriptstyle{
    x_{0} x_{45} + x_{1} x_{26} - x_{6} x_{20} + x_{8} x_{18} -     x_{10} x_{16} + x_{12} x_{13}=0, }\\   \scriptstyle{
    x_{1} x_{48} + x_{2} x_{50} + x_{3} x_{52} + x_{4} x_{53} +    x_{6} x_{54} + x_{7} x_{55}=0, }&  \scriptstyle{
    x_{0} x_{46} - x_{2} x_{24} + x_{3} x_{22} - x_{5} x_{19} +     x_{8} x_{17} - x_{9} x_{15}=0, }\\   \scriptstyle{
    x_{0} x_{47} - x_{1} x_{25} + x_{4} x_{20} - x_{5} x_{18} +      x_{10} x_{14} - x_{11} x_{12}=0, }&  \scriptstyle{
    x_{0} x_{48} + x_{2} x_{23} - x_{3} x_{21} + x_{4} x_{19} -      x_{6} x_{17} + x_{7} x_{15}=0, }\\   \scriptstyle{
    x_{0} x_{49} + x_{1} x_{24} - x_{3} x_{20} + x_{5} x_{16} -      x_{8} x_{14} + x_{9} x_{12}=0, }&  \scriptstyle{
    x_{0} x_{50} - x_{1} x_{23} + x_{3} x_{18} - x_{4} x_{16} +      x_{6} x_{14} - x_{7} x_{12}=0, }\\   \scriptstyle{
    x_{0} x_{51} - x_{1} x_{22} + x_{2} x_{20} - x_{5} x_{13} +      x_{8} x_{11} - x_{9} x_{10}=0, }&  \scriptstyle{
    x_{0} x_{52} + x_{1} x_{21} - x_{2} x_{18} + x_{4} x_{13} -      x_{6} x_{11} + x_{7} x_{10}=0, }\\   \scriptstyle{
    x_{0} x_{53} - x_{1} x_{19} + x_{2} x_{16} - x_{3} x_{13} +      x_{6} x_{9} - x_{7} x_{8}=0, }&  \scriptstyle{
    x_{0} x_{54} + x_{1} x_{17} - x_{2} x_{14} + x_{3} x_{11} -      x_{4} x_{9} + x_{5} x_{7}=0, }\\   \scriptstyle{
    x_{0} x_{55} - x_{1} x_{15} + x_{2} x_{12} - x_{3} x_{10} +      x_{4} x_{8} - x_{5} x_{6}=0, }
\end{array}
$$
\begin{tabular}{p{1\textwidth}}
 $\scriptstyle{x_{0} x_{28} + x_{1} x_{29} + x_{2} x_{30} + x_{3} x_{31} + 
        x_{4} x_{32} + x_{5} x_{33} + x_{6} x_{34} + x_{8} x_{36} + 
        x_{10} x_{38} + x_{12} x_{40} + x_{15} x_{43} - 
        x_{27} x_{55}=0, } $\\$  \scriptstyle{
    3 x_{0} x_{28} + 3 x_{1} x_{29} + 3 x_{2} x_{30} + 3 x_{3} x_{31}
        + 3 x_{4} x_{32} + x_{5} x_{33} + 3 x_{6} x_{34} + 
        3 x_{7} x_{35} + x_{8} x_{36} + x_{9} x_{37} + x_{10} x_{38} 
        + x_{11} x_{39} +}$ $\scriptstyle{
  + x_{12} x_{40} + x_{13} x_{41} + 
        x_{14} x_{42} + x_{15} x_{43} + x_{16} x_{44} + x_{17} x_{45}
        + x_{18} x_{46} + x_{19} x_{47} - x_{20} x_{48} + 
        x_{21} x_{49} - x_{22} x_{50} + x_{23} x_{51}-}$ $\scriptstyle{
- x_{24} x_{52}
        - x_{25} x_{53} - x_{26} x_{54} - x_{27} x_{55}=0, }$\\$  \scriptstyle{
    2 x_{0} x_{28} + 2 x_{1} x_{29} + 2 x_{2} x_{30} + 2 x_{3} x_{31}
        + 2 x_{4} x_{32} + 2 x_{5} x_{33} + x_{6} x_{34} + 
        x_{7} x_{35} + x_{8} x_{36} + x_{9} x_{37} + x_{10} x_{38} + 
        x_{11} x_{39} +}$ $\scriptstyle{
+ x_{12} x_{40} + x_{14} x_{42} + x_{15} x_{43}
        + x_{17} x_{45} - x_{26} x_{54} - x_{27} x_{55}=0, }$\\$  \scriptstyle{
    3 x_{0} x_{28} + 3 x_{1} x_{29} + 3 x_{2} x_{30} + 3 x_{3} x_{31}
        + 2 x_{4} x_{32} + 2 x_{5} x_{33} + 2 x_{6} x_{34} + 
        2 x_{7} x_{35} + 2 x_{8} x_{36} + 2 x_{9} x_{37} + 
        x_{10} x_{38} + x_{11} x_{39} +}$ $\scriptstyle{
+ x_{12} x_{40} + x_{13} x_{41}
        + x_{14} x_{42} + x_{15} x_{43} + x_{16} x_{44} + 
        x_{17} x_{45} + x_{19} x_{47} - x_{25} x_{53} - x_{26} x_{54}
        - x_{27} x_{55}=0, }$\\$  \scriptstyle{
    5 x_{0} x_{28} + 5 x_{1} x_{29} + 5 x_{2} x_{30} + 3 x_{3} x_{31}
        + 3 x_{4} x_{32} + 3 x_{5} x_{33} + 3 x_{6} x_{34} + 
        3 x_{7} x_{35} + 3 x_{8} x_{36} + 3 x_{9} x_{37} + 
        3 x_{10} x_{38} + 3 x_{11} x_{39} +}$ $\scriptstyle{
+ x_{12} x_{40} + 
        3 x_{13} x_{41} + x_{14} x_{42} + x_{15} x_{43} + 
        x_{16} x_{44} + x_{17} x_{45} + x_{18} x_{46} + x_{19} x_{47}
        + x_{20} x_{48} + x_{21} x_{49} + x_{22} x_{50} - 
        x_{23} x_{51}-}$ $\scriptstyle{
 - x_{24} x_{52} - x_{25} x_{53} - x_{26} x_{54}
        - x_{27} x_{55}=0, }$\\$  \scriptstyle{
    2 x_{0} x_{28} + 2 x_{1} x_{29} + x_{2} x_{30} + x_{3} x_{31} + 
        x_{4} x_{32} + x_{5} x_{33} + x_{6} x_{34} + x_{7} x_{35} + 
        x_{8} x_{36} + x_{9} x_{37} + x_{10} x_{38} + x_{11} x_{39} + }$ $\scriptstyle{
 +       x_{12} x_{40} + x_{13} x_{41} + x_{14} x_{42} + x_{16} x_{44}
        + x_{18} x_{46} + x_{20} x_{48}=0, }$\\$  \scriptstyle{
    3 x_{0} x_{28} + x_{1} x_{29} + x_{2} x_{30} + x_{3} x_{31} + 
        x_{4} x_{32} + x_{5} x_{33} + x_{6} x_{34} + x_{7} x_{35} + 
        x_{8} x_{36} + x_{9} x_{37} + x_{10} x_{38} + x_{11} x_{39} +}$ $\scriptstyle{
        +x_{12} x_{40} + x_{13} x_{41} + x_{14} x_{42} + x_{15} x_{43}
        + x_{16} x_{44} + x_{17} x_{45} + x_{18} x_{46} + 
        x_{19} x_{47} + x_{20} x_{48} + x_{21} x_{49} + x_{22} x_{50}
        + x_{23} x_{51} +}$ $\scriptstyle{
       + x_{24} x_{52} + x_{25} x_{53} + 
        x_{26} x_{54} + x_{27} x_{55}}$
\end{tabular}

Again the equations encode the $\gote_7$ Lie Algebra.

\subsection{Degenerate cases}
In all the above examples the Lie algebra $\gotg$ obtained from the quadrics in the ideal of $X$ is exactly the maximal
 semisimple subalgebra of the parabolic subalgebra of the original nilpotent orbit. This happens for the two degenerate
 cases as well.\\
\begin{figure}[htb]
\centering
\includegraphics[width=0.9\textwidth]{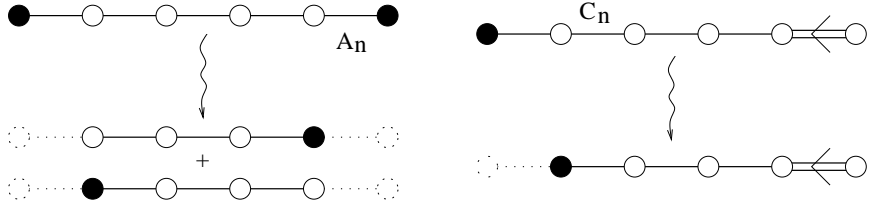}
\caption{The degenerate cases: $A_{n+1}$ gives $A_{n-1} \cup A_{n-1}$ and $C_{n+1}$ gives $C_{n}$.}
 \label{pic_adjoint_degenerate}
\end{figure}

For $\gotsl_{n+2}$ the contact manifold is $\P(T^*\P^{n+1})$ and the legendrian variety is a disjoint union of two
linear $\P^{n-1}$'s embedded in $\P^{2n-1}$. Therefore choosing proper symplectic coordinates the ideal is generated
by $\left\{x_i x_{n+j} : i,j \in \{0,1,\ldots n-1\}\right\}$ and this gives $\gotg \simeq \gotsl_n$ which is exactly the 
maximal semisimple subalgebra of $\gotp < \gotsl_{n+2}$.\\

For $\gotsp_{2n+2}$ the contact manifold is $\P^{2n+1}$ and the corresponding ``legendrian variety'' is empty.
So the quadratic part of the ideal of a empty variety is all the $\Sym^2V^*$ and hence $\gotg \simeq \gotsp_{2n}$.
This is again the maximal semisimple subalgebra in $\gotp < \gotsp_{2n+2}$.


\section{Lie groups and the representation theory}\label{representation_theory}

I will summarize several facts about the Lie groups and the representation theory.
They are  essential to the final classification in the preceding section. 
Most of them is well known, but I was not able to find the appropriate references to such explicit statements.\\

\subsection{Semisimple and irreducible}

The first statement  says about the semisimplicity of the group, 
which was defined in the section \ref{rozdzial_rozm_gen_kwadryki}.

\begin{lem}\label{semisimple}
 Let $X$ be a projective variety and let $G$ be a group of projective automorphisms preserving $X$. Suppose that $G$ acts
transitively on $X$ and that the subgroup acting trivially on $X$ is discrete. Then $G$ is semisimple.
\end{lem}
\begin{prf}
Since $X$ is projective, there exist a parabolic subgroup $P<G$ (see \cite[\S11.2]{borel}), such that
$G \slash P = X$. But since the subgroup acting
trivially is discrete, then so is the intersection  $ \bigcap_{g \in G} g^{-1} P g$. Now $P$ contains a Borel subgroup $B$
(\cite[cor. 11.2]{borel}),
so $\bigcap_{g \in G} g^{-1} B g$ is discrete as well and therefore it's connected component of the identity
(which is equal to the radical of $G$ - see \cite[\S11.21 and thm 11.1]{borel}) is
trivial, so $G $ is semisimple (see \cite[\S11.21]{borel}).
\end{prf}

Recall, that every representation of a semisimple Lie algebra decomposes into eigenspaces of the action of its Cartan
subalgebra $\goth < \gotg$: 
$$
\gotg = \goth \oplus \bigoplus_{\alpha \in R} \gotg_\alpha
 \ \ \ \
V =  \bigoplus_{\lam \in \goth^*} V_\lam
$$
where $R$ is the set of roots and $\gotg_\alpha$ is the root space of $\alpha$. 
For a convenience I will also denote by $\gotg_0$ the 0-eigenspace,
i.e. the Cartan subalgebra $\goth$, although $0$ is not considered as a root.
V is a representation and $V_\lam$ is the $\lam$-eigenspace.

Next statement says about the irreducibility of the restriction of the natural representation
 of $\gotsp(V)$ to the Lie algebra $\gotg$.

\begin{lem}\label{irreducible}
Let $G$ be a semisimple group, $V$ its representation and $X\subset \P(V)$ a closed orbit of the induced action on $\P(V)$.
 Then the linear space $W:=\sspan \hat{X}$ is an irreducible subrepresentation of $V$. 
\end{lem}
\begin{prf}
Clearly $\sspan \hat{X}$ is a subrepresentation.
Choose a point $x \in X$ and let $P<G$ be the parabolic subgroup preserving $x$.
 Let also $\gotp$ and $\gotg$ be the
corresponding Lie algebras. Also choose a Cartan subalgebra $\goth < \gotp < \gotg$ \cite[\S23.3]{fultonharris}
 and a suitable root order.
Then all the elements in any positive root space of $\gotg$ take $x$ to $0$ 
and therefore a nonzero vector in the line $x \subset V$
is the highest weight vector for a irreducible subrepresentation $W$. But then $x \in \P(W)$ and hence $W$ must contain
$\sspan \hat{X}$. So $W=\sspan \hat{X}$ and it is a irreducible subrepresentation.
\end{prf}

\subsection{Preferred basis}

Assume that $G$ is a semisimple group, $V$ is an irreducible representation and $X$ is the closed orbit in $\P(V)$.
Let $n-1$ be the dimension of $X$ and $m$ be the dimension of $V$. 
Fix a parabolic subgroup $P<G$ such that $X \simeq G \slash P$, a Cartan subalgebra $\goth < \gotp < \gotg$ and 
the corresponding root order as in the above proof. Then:

\begin{lem}\label{basis_of_representation}
There exist a basis $\{v_0, v_1,\ldots, v_{n-1}, v_n, \ldots, v_{m-1}\}$ of $V$
\begin{itemize}
\item[1)] the action of $\goth$ on $V$ diagonalises in the basis;
\item[2)] $v_0$ is the highest weight vector and $v_0\in \hat{X}$; 
\item[3)] the tangent space of $\hat{X}$ at the point $v_0$ is spanned by $v_0, v_1, \ldots, v_{n-1}$.
\end{itemize}
\end{lem}
\begin{prf}
Choose a basis satisfying 1) and 2). Then the tangent space to $\hat{X}$ at $v_0$ is exactly the image $\gotg(v_0)$.
So let $\alpha_1,\ldots, \alpha_{n-1}$ be the negative roots that are not is $\gotp$ and let $g_1, \ldots, g_{n-1}$
denote the corresponding Lie algebra elements. Then $v_0$ and $g_i(v_0)$'s are weight vectors, which span the tangent 
space. So let $v_i:=g_i(v_0)$ and choose $v_n, \ldots, v_{m-1}$ out of the other weight vectors correspondingly. 
\end{prf}

\subsection{Semisimple, but not simple}
Further, I analyze the case where $G$ is semisimple but not simple. Then a nice statement can be formed 
about the irreducible representations of such.

\begin{lem}\label{lemma_for_semisimple}
Let $\gotg=\gota \oplus \gotb$ for some semisimple Lie algebras $\gota, \gotb$ 
and let $V$ be an irreducible representation of $\gotg$
Then $V \simeq W_\gota \otimes W_\gotb$, where $W_\gota$ is an irreducible representation of $\gota$ and $W_\gotb$ is an
 irreducible representation of $\gotb$.
\end{lem}
\begin{prf}
First decompose the highest weight $\lam_0$ of $V$ into the components: $$\lam_0 = \alpha_0 + \beta_0,$$
$\alpha_0$ being a weight of $\gota$ and $\beta_0$ being a weight of $\gotb$. Next let $W_\gota$ be the irreducible
representation with the highest weight $\alpha_0$ and $W_\gotb$ be the irreducible
representation with the highest weight $\beta_0$ and denote the highest weight vectors of these representations
via $w_{\alpha_0}$ and $w_{\beta_0}$ correspondingly.
 Clearly $W:=W_\gota \otimes W_\gotb$ contains the representation $V$ as a
subrepresentation, so the only thing is to prove that $W$ is irreducible. Now $W$ is spanned by simple tensors of weight
vectors, i.e. by the tensors of the form $w_\alpha \otimes w_\beta$ for $w_\alpha$ a weight vector in $W_\gota$ and
$w_\beta$ a weight vector in $W_\gotb$. But since $W_\gota$ (and $W_\gotb$) is irreducible, it follows that
$w_\alpha= a_k \circ \ldots \circ a_1 (w_{\alpha_0})$ for some $a_i$'s in some root spaces of $\gota$ (and similarly
$w_\beta= b_l \circ \ldots \circ b_1 (w_{\beta_0})$ for some $b_j$'s in some root spaces of $\gotb$).
Then
$$w_\alpha \otimes w_\beta= b_l \circ \ldots \circ b_1 \circ a_k \circ \ldots \circ a_1 (w_{\alpha_0} \otimes w_{\beta_0})
$$
and hence $W$ is irreducible, so $V \simeq W= W_\gota \otimes W_\gotb$.
\end{prf}

\subsection{Semisimple subalgebra generated by $\alpha_1, \ldots, \alpha_k$}\label{subalgebra_s}

Assume again that $\gotg$ is a semisimple Lie algebra. Recall that for every (positive) root $\alpha \in R$ one 
has a subalgebra $\gots_\alpha < \gotg$ spanned by the root spaces $\gotg_\alpha$, $\gotg_{-\alpha}$ and their
product $[\gotg_\alpha,\gotg_{-\alpha}]$. Every subalgebra $\gots_\alpha$ is isomorphic to $\gotsl_2$ 
\cite[fact 14.6 and (D.16)]{fultonharris}
 and a lot of interesting properties of representations of semisimple Lie algebras can be proved just by restricting to 
$\gots_\alpha$'s. For example:

\begin{lem}\label{g_map_non_zero}
Let $V$ be a representation of $\gotg$, $\alpha$ be a root of $\gotg$ and $\lam$ be a weight of $V$. 
Assume that $\lam + \alpha$ is a weight of $V$ as well.
Let $g$ be a non-zero element of the root space $\gotg_\alpha$. Then the linear map
$$g|_{V_\lam} : V_\lam \lra V_{\lam + \alpha} $$
is non-zero.
\end{lem}

\begin{prf}
Let $W= \bigoplus_{t \in \R} V_{\lam + t \alpha}$, i.e. the sum of those weight spaces, which are in the 
``$\alpha$ string through $\lam$'' 
(i.e. the intersection of weights of $V$ and the line passing through $\lam$ and parallel to 
$\alpha$ - see for example \cite[\S21.1, property (5)]{fultonharris} 
for an analogous notion for the adjoint representation).
Then $W$ is a representation of $\gots_\alpha$ and the problem is reduced to the case of $\gotg \simeq \gotsl_2$.
And for $\gotsl_2$ it is easy - see \cite[(11.5)]{fultonharris}.
\end{prf}

Now I will present an obvious generalization of the above idea. So let $\alpha_1, \ldots, \alpha_k \in R$ be some
positive and linearly independent roots of $\gotg$ and let $S\subset R$
be the intersection of $R$ and $S_{\R} := \sspan_\R \{\alpha_1, \ldots, \alpha_k\}$.
Let $\gots=\gots_{\alpha_1, \ldots, \alpha_k}$
 be the subalgebra of $\gotg$ generated by the root spaces corresponding to the roots of $S$. 
Then $\gots$ is a semisimple Lie algebra of rank k.

For short, I will refer to the subalgebra $\gots$ as to the \textbf{semisimple subalgebra generated by
$\alpha_1, \ldots, \alpha_k$}, although it is not perfectly precise.
Note that for $k=1$ the subalgebra $\gots_{\alpha_1}$ is nothing else than the
$\gotsl_2$ subalgebra corresponding to the root $\alpha_1$.

Further let $V$ be a representation of $\gotg$. If $V_{\gots}$ is the  restriction of  the representation to $\gots$, 
it splits into the \textbf{slice-subrepresentations} (usually not irreducible),
 each of which has the weights (and multiplicities) being a slice of the weight
diagram  of $V$, i.e. the intersection of the weight diagram and $\lam+ S_{\R}$ ($\lam$ is a weight of $V$).
Now in particular:
\begin{lem}
 If $V$ is irreducible take the highest weight $\lam_0$ and let $W$ be the irreducible representation
 of $\gots$ with the highest weight $\lam_0|_{\gots}$. Then the weights of $W$ are contained in those weights of $V$
 that belong to $\lam_0 + S_{\R}$ and the multiplicities of $W$ are less or equal than those of the slice.
\end{lem}
\noprf

\subsection{Nilpotent subalgebra}

Let $\gotn<\gotg$ be the nilpotent subalgebra generated by all root spaces $\gotg_{\alpha_i}$
where the $\alpha_i$'s are these roots which are not the roots of $\gotp$. 
So $\gotg= \gotp \oplus \gotn$ (the sum is of the vector spaces, not of the Lie Algebras)
 and $R= R_\gotp \cup R_\gotn$, where $R_\gotp$ and $R_\gotn$ are
the sets of roots of $\gotp$ and $\gotn$ correspondingly.

Let $D$ be the Dynkin diagram for $\gotg$.
Recall, that a homogeneous space is determined by specifying which of the simple negative roots belong to
$R_\gotn$ - all the others belong to $R_\gotp$. 
\begin{prop}\label{one_simple_root_generally}
Assume $\gotg$ is simple. Also assume that the angle between any two roots in $R_\gotn$ is either acute or right. 
Then only there is only one simple root in $R_\gotn$.
\end{prop}
\begin{prf}
Assume that there are two simple roots $\alpha$ and $\alpha'$ in $R_\gotn$.
Since $\gotg$ is simple, the Dynkin diagram $D$ is connected.
So let
$$\alpha=\alpha_0, \alpha_1, \ldots, \alpha_{k-1}, \alpha_k = \alpha'$$
be a connected string of negative simple roots (without repetitions).
\begin{lem}
For every $0\le m < k$ the sum $\alpha_0+\ldots+\alpha_m$ is a root in $R_\gotn$.
\end{lem}
\begin{prf}
Argue inductively.
For $m=0$ there is nothing to do. So assume $\alpha_0+\ldots+\alpha_{m-1} =: \beta$ is in $R_\gotn$ 
for some $m \ge 1$. 
Then the angle between $\beta$ and $\alpha_m$ is obtuse and hence  $\beta + \alpha_m$ is a root 
\cite[\S21.1, property (6), p324]{fultonharris}. 
Also it has to  be in $R_\gotn$, for otherwise:
$$\beta = ( \beta + \alpha_m) + (-\alpha_m) $$
is a sum of two roots in $R_\gotp$, so $\beta$ is in $R_\gotp$, contradicting the inductive assumption.
This proves of the lemma.
\end{prf}
To finish the proof of the proposition just notice that the angle between $\alpha_0+\ldots+\alpha_{k-1}$ and $\alpha_k$ is
obtuse, as in the proof of the lemma and that contradicts the assumption of the proposition. 
\end{prf}


\section{Legendrian subvarieties generated by quadrics}\label{quadrics}

The theorems of section \ref{rozdzial_rozm_gen_kwadryki} give a powerful tool for explicit calculations 
for legendrian varieties generated by quadrics. Namely, I can study the representation theory for the group $G$.
Using the representation theory I can finally restrict the number of possibilities just to the cases described in 
the section \ref{rozdzial_przyklad}.

\subsection{Notation and basic properties}
In this subsection I want to fix a notation for whole section \ref{quadrics}. So let $X\subset_l \P(V)$ be a smooth 
legendrian variety generated by quadrics (also I want $X$ to be irreducible and nondegenerated). Recall from 
the section \ref{rozdzial_rozm_gen_kwadryki} the subgroup $G< \Sp$, which acts transitively on $X$ by the theorem 
\ref{rozmaitosci_jednorodne}. Let $\gotg$ be it's Lie algebra and $\I_2$ the set of quadrics in the ideal of $X$.

\begin{cor}\label{G_is_semisimple}
With the above notation, $G$ is semisimple and $V$ is its irreducible representation.
\end{cor}

\begin{prf}
It is an obvious consequence of the above assumptions,  the proposition \ref{kernel_of_action} and finally the 
lemmas \ref{semisimple}  and \ref{irreducible}. 
\end{prf}

This is already a lot. The forthcoming corollary sketches the idea of the final classification in the subsection 
\ref{section_final_classification}. So a list of properties of the representation $V$ will be proved.
Further, I will check that only the examples of section \ref{rozdzial_przyklad} satisfy the properties. 
The irreducibility is more than enough to deal with the $G_2$ case.

\begin{cor}\label{not_g2}
$G$ is not of type $G_2$.
\end{cor}
\begin{prf}
Assume $\gotg$ is isomorphic to the exceptional Lie algebra $\gotg_2$ (\cite[lecture 22]{fultonharris},
\cite[\S19.3]{humphreys}).
There are exactly three homogeneous spaces for the group $G_2$: the five dimensional quadratic hypersurface,
the five dimensional $G_2$-variety and the six dimensional flag over the two previous spaces - i.e. the quotient by the
Borel subgroup (\cite[\S23.3, p.391]{fultonharris}). Therefore, since dimension of
$V$ is twice the dimension of $\hat{X}$, it follows that $\dim V$ is either $12$ or $14$. Have a quick look at
the representation theory of $\gotg_2$  (\cite[\S22.3]{fultonharris})
and notice that there are only two non-trivial irreducible representations of dimension
less or equal to 14: the natural one (of dimension 7) and the adjoint one (of dimension 14). So $V$ must be the adjoint
representation of $\gotg_2$, but then $X$ would have to be the $G_2$-variety
(\cite[\S23.3, p.391]{fultonharris}), so the dimension argument fails. Hence
$\gotg$ is not $\gotg_2$ nor $G$ is of type $G_2$.
\end{prf}

Note, that it is not essential to put the above corollary here, yet it will  simplify some proofs of the subsection
\ref{section_final_classification}.

\begin{prop}\label{basis_of_V}
There exist a basis $\{v_0, v_1,\ldots, v_{n-1}, v_n, \ldots, v_{2n-1}\}$ of $V$ and a Cartan subalgebra  $\goth < \gotg$,
 such that:
\begin{itemize}
\item[1)] the action of $\goth$ on $V$ diagonalises in the basis;
\item[2)] $v_0 \in \hat{X}$  and $v_0$ is the highest weight vector of $V$;
\item[3)] the tangent space of $\hat{X}$ at the point $v_0$ is spanned by $v_0, v_1, \ldots, v_{n-1}$.
\item[4)] the basis is symplectic, i.e. the matrix of the symplectic form is:
$$J=\begin{pmatrix}
  0 & \Id_n \\
  -\Id_n & 0
\end{pmatrix}.$$
\end{itemize}
\end{prop}
\begin{prf}
Let $\{v_0,\ldots, v_{n-1}, v_{n}, \ldots,  v_{2n-1}\}$ be a basis of $V$ and $\goth$ a Cartan subalgebra satisfying 
1), 2) and 3) - see the lemma \ref{basis_of_representation}. In fact, this basis is already
little away from being symplectic. To see that, for $h\in \goth$ write
$$h= \diag(\lam_0(h),\lam_1(h),\ldots, \lam_{2n-1}(h)) $$
where $\lam_i$'s are linear forms on $\goth$. One has the following equalities for every \mbox{$0\le i,j \le 2n-1 $}:
$$
\lam_i(h) \cdot \omega(w_i,w_j)= \omega(h (w_i),w_j)= - \omega(w_i,h (w_j) )= -\lam_j(h) \cdot \omega(w_i,w_j).
$$
Therefore, whenever $\omega (w_i,w_j)$ is non-zero, one has $\lam_i = -\lam_j$. Using this property and the definition of
the symplectic form one can easily change the order of $v_i$'s and rescale them 
(and also perform another base change within weight spaces)
 to get the property 4) of the proposition without spoiling 1), 2) and 3). 
\end{prf}

From now on I will stick to the \textbf{notation}, that $\goth$ is the Cartan algebra and   
$\{v_0, v_1,\ldots, v_{n-1}, v_n, \ldots v_{2n-1}\}$ is the 
basis of the proposition \ref{basis_of_V}. Note that still we have a small choice of this basis: namely we can permute 
$v_1,\ldots v_{n-1}$ and apply the same permutation to  $v_{n+1},\ldots, v_{2n-1}$ and the properties 1)-4) of the 
proposition will be preserved.

Now, thinking of $\gotg$ as a subalgebra of $\gotsp_{2n}$, take any $g \in \gotg$ and write it as a matrix in the basis 
$$\{v_0, v_1,\ldots, v_{n-1}, v_n, \ldots v_{2n-1}\}$$ in the following block form: 

\begin{equation} \label{form_of_matrix_g_}
\left( \begin{array}{c|ccc|c|ccc}
\lam_0 && a_2^T && \nu     && c^T&    \\
\hline
       &&       &&         &&    &    \\
a_1    && A     && c       && C  &    \\
       &&       &&         &&    &    \\
\hline
\mu    && b^T   && -\lam_0 && -a_1^T& \\
\hline
       &&       &&         &&    &    \\
b      && B     && -a_2    && -A^T  &     \\
       &&       &&         &&    &    \\
\end{array} \right)
\end{equation}
where $\lam_0, \mu$ and $\nu$ are scalars of $\F$, $a_1,a_2, b$ and $c$ are vertical vectors in $\F^{n-1}$ and 
finally $A,B$ and $C$ are $(n-1) \times (n-1)$ matrices. Moreover $B$ and  $C$ are symmetric matrices. For a convenience,
I will write for example $\lam_0(g)$ or $c(g)$ or $C(g)\ldots$ to mean the proper block of the matrix $g$. 

Let me state the first properties of the above block form of $g$:

\begin{lem}\label{properties_of_g}
With the above setup, the following conditions hold:
\begin{itemize}
\item[(i)]  $\mu=0$ and $b = 0$;
\item[(ii)] $\lam_0$ and $a_1$ are epimorphic, i.e. for each $ \lam \in \F$ and $a \in \F^{n-1}$ there exists 
$g \in \gotg$, such that $\lam_0(g)=\lam$ and $a_1(g) = a$.
\end{itemize}
\end{lem}

\begin{prf}
To prove (i) first notice, that since $v_0\in \hat{X}$ and the quadric corresponding to $2 J \circ g$ is in $\I_2$ (see 
\ref{isomorphism_Sym2_sp}), hence it 
vanishes on $v_0$ and this is equivalent to $\mu=0$.
 
Second, to prove the other part of (i) and (ii) as well notice, that by computing the linear part of $\exp(g)$, 
we have that $T_{v_0}\hat{X}$ is exactly equal to $$ \{ (\lam_0(g), a_1(g), \mu (g), b(g)) : g \in \gotg \} . $$
 On the other hand by the property 4) of the basis, $$T_{v_0}\hat{X} =\sspan \{ v_0, v_1, \ldots, v_{n-1}\} . $$
So $b=0$ and both $\lam_0$ and $a_1$ can be chosen arbitrarily. 
\end{prf}

Summarizing we get the following block form of the matrix $g$:

\begin{equation}  \label{form_of_matrix_g}
\left( \begin{array}{c|ccc|c|ccc}
\lam_0(g) && a_2(g)^T && \nu(g)  && c(g)^T&    \\
\hline
          &&            &&          &&       &    \\
a_1(g)    && A(g)       && c(g)     && C(g)  &    \\
          &&            &&          &&       &    \\
\hline
 0        && 0          &&-\lam_0(g)&&-a_1(g)^T&  \\
\hline
          &&            &&          &&       &    \\
0         &&      B(g)  && -a_2(g)  &&-A(g)^T&    \\
          &&            &&          &&       &    \\
\end{array} \right) . 
\end{equation}

Next, for an element $h \in \goth$ write  
$$h= \diag (\lam_0(h),\lam_1(h),...\lam_{n-1}(h), -\lam_0(h),-\lam_1(h),...-\lam_{n-1}(h))$$ 
for some $\lam_i \in \goth^*, i \in \{0,1,\ldots n-1 \} $. Note, that obviously $\lam_0$ coincides with
the restriction of the form defined on $\gotg$ in \ref{form_of_matrix_g_}.

\subsection{Roots and weights}

In this subsection I am going to examine the properties of roots of $\gotg$ and weights of $V$.
So denote by $R$ the set of roots with respect to the Cartan subalgebra $\goth$.

\begin{prop}\label{weights_and_root_of_V}
\begin{itemize}
\item[(1)] The weights of $V$ are exactly $\pm \lam_i$'s  for $i \in \{0,\ldots n-1 \}$. 
In particular multiple weights occur if and only if $\lam_i = \pm \lam_j$ for some $i,j$ and some choice of sign.
\item[(2)] Every root of $\gotg$ is of the form $\pm \lam_i \pm \lam_j$ for some choice of $i$ and $j$ and 
some choice of signs - but the choices are never unique.
\item[(3)] $\lam_0$ is the highest weight of $V$.
\end{itemize}
\end{prop}

\begin{prf}
 (1) is obvious and (3) follows immediately from  the point 2) of the proposition \ref{basis_of_V}. 
To prove the first part of (2)
notice that the adjoint representation of $\gotg$ is a subrepresentation of the adjoint representation of 
$\gotsp(V)$ restricted to $\gotg$. The last one is isomorphic to $\Sym^2 V^* \simeq \Sym^2 V$ and the weights of
$\Sym^2 V$ are exactly $\pm \lam_i \pm \lam_j$.

The first part of (2) could also be seen explicitly by computing the matrix form of $[h,g]$ for $h \in \goth$, 
$g \in \gotg$.

 The second part of (2) is also easy: if there is a root $\alpha$ with unique presentation $\alpha=\pm \lam_i \pm \lam_j$
then the matrix of every $g \in \gotg_\alpha$ is of rank at most 2. But matrices of rank 2 correspond to reducible 
quadrics in the ideal of $X$, contradicting the assumptions of irreducibility or nondegeneracy.
\end{prf}

\begin{cor}\label{vn_in_X}
$2\lam_0$ is not a root of $\gotg$, $\nu=0$ and $v_n \in \hat{X}$.
\end{cor}
\begin{prf}
Since $\lam_0$ is the highest weight of $V$,
$2\lam_0$ is not equal to any other sum of the form $\pm \lam_i \pm \lam_j$. So it cannot be a root of $\gotg$.

Next, this means, that for every $g\in \gotg_\alpha$ (for every $\alpha$) $\nu(g)=0$, so simply $\nu =0$. Computing
the quadric corresponding to $Jg$, one
can see that $\nu=0$ is equivalent to the fact that all the quadrics in $\I_2$ vanish on $v_n$. So $v_n\in \hat{X}$.
\end{prf}

\begin{lem}\label{lambdas_non_zero}
All $\lam_i$'s are non-zero forms on $\goth$.
\end{lem}
\begin{prf}
Assume $\lam_i=0$. For any $g\in \gotg$ let $g_{(j,k)}$ be the term in the $j$'th row and $k$'th column of $g$
(for a convenience I enumerate rows and columns from 0). Then:
\begin{equation}\label{equation_for_lambdas_non_zero}
 g_{(i,i)}= g_{(i,n+i)} = g_{(n+i,i)} =g_{(n+i,n+i)} =0.
\end{equation}
Indeed, if any of the above terms is non-zero for some $g \in \gotg$, then it is also non-zero for some 
$g_\alpha \in \gotg_\alpha$,
for some $\alpha \in R \cup \{0\}$. But then computing the action of $\goth$ on $g_\alpha$ and using the assumption that
$\lam_i = 0$, I get that $\alpha =0$. So in fact $g\in \goth$ and this contradicts either the assumption that $\lam_i=0$ or
that every element in $\goth$ is a diagonal matrix.

But translating the equations \eqref{equation_for_lambdas_non_zero} into the properties of the ideal $\I$ of $X$
one gets that the line spanned by $\{v_i, v_{n+i} \}$ is contained in $X$ contradicting the basic properties
of  legendrian variety.
\end{prf}

Now,  choose an order of the roots $R=R^+ \cup R^-$. It must be compatibile with the setup so far,
 so that the parabolic subalgebra $\gotp$ fixing $[v_0] \in \P(V)$ must contain the Borel subalgebra 
$\goth \oplus \bigoplus_{\alpha \in R^+} \gotg_\alpha \subset \gotp$.

Note, that if $\alpha$ is a root  then every element $g\in \gotg_\alpha$ 
has only $0$ terms on the diagonal. Similarly if the term of $g$ in $i$'th 
row and $j$'th column (again for a convenience I enumerate rows and columns from 0)
is non zero, then:
$$ \begin{array}{cl}
\alpha = \lam_i - \lam_j & \textrm{for } 0 \le i,j \le n-1, \\
\alpha = \lam_i + \lam_j & \textrm{for } 0 \le i \le n-1, \ n \le j \le 2n-1 \textrm{ and} \\
\alpha =-\lam_i - \lam_j & \textrm{for } n \le i \le 2n-1, \ 0 \le j \le n-1.
\end{array}
$$
In particular taking $j=0$ and $1 \le i \le n-1$ and using the lemma \ref{properties_of_g}(ii) we get that:
\begin{prop}
Each $\lam_i - \lam_0$ is a root of $\gotg$. Moreover $\lam_i - \lam_0 \in R^-$ and each of these roots is 
different.
\end{prop}
\begin{prf}
To see that  $\lam_i - \lam_0 \in R^-$, just notice, that if $a_1(g) \ne 0$ then $g \notin \gotp$. To see that for 
$i\ne j$ the roots $\lam_i - \lam_0$ and $\lam_j - \lam_0$ are different, it suffices to use the fact that the root spaces
are one dimensional (see \cite[fact 14.2(i) and (D.20)]{fultonharris} or  \cite[prop. 8.4(a)]{humphreys}).  
\end{prf}

\begin{cor} \label{lambdas_are_different}
For every $i,j \in \{0,1,\ldots,n-1\}$, $i\ne j$ the weights $\lam_i$ and $\lam_j$ are different. 
\end{cor}
\noprf

\begin{cor}\label{positive_roots}
For each $1 \le i \le n-1$ the form $\lam_0 - \lam_i$ is a positive root.
\end{cor}
\begin{prf}
The fact that the negative of a root is a root can be found for example in \cite[fact 14.2(iii) and (D.13)]{fultonharris}
 or  \cite[thm. 8.5(b)]{humphreys}.
\end{prf}

Now let me study in more details the structure of $\gotg_{\lam_0-\lam_i}$. According to what has been done so far every 
element $g_i\in \gotg_{\lam_0-\lam_i}$ is of the following form:
$$
\left( \begin{array}{c|ccc|c|ccc}
0 && k_i(g_i)e_i^T&&     0                    &&c(g_i)^T &  \\
\hline 
  &&                            &&                          &&                       &  \\
0 && A(g_i)       &&  c(g_i)    &&C(g_i)   &  \\
  &&                            &&                          &&                       &  \\
\hline
0 &&             0              &&   0                      &&                      0&  \\
\hline
  &&                            &&                          &&                       &  \\
0 &&     B(g_i)                 &&-k_i(g_i)e_i&&-A(g_i)^T&    \\
  &&                            &&                          &&                       &    \\
\end{array} \right) 
$$
where $k_i(g_i)e_i$ is just a scalar multiple of $e_i$ and $e_i$ is a $i^{\textrm{th}}$ standard base vector 
of $\C^{n-1}$. Just in case the reader feels a bit lost at this moment, let me quickly recall, how did I get the above 
form out of \eqref{form_of_matrix_g}:
\begin{itemize}
\item 
$\lam_0(g_i)=0$, because $\lam_0-\lam_i$ is not 0, so $g_i \notin \goth$;
\item
$a_1(g_i)=0$, because $\lam_0-\lam_i \ne \lam_j - \lam_0$ -- the first one is positive, while the other is a
negative root;
\item
$a_2(g_{i})$ is a multiple of $e_i$, because $\lam_0-\lam_i \ne \lam_0 - \lam_j$ for $j\ne i$ -- that is the
 corollary \ref{lambdas_are_different}; now the goal is to prove that $k_i$ actually is non-zero;
\item
$\nu(g_i)=0$, because $\nu=0$ -- that is the corollary \ref{vn_in_X}.  
\end{itemize}

Now assume $g_i$ is non-zero. Then by the lemma \ref{g_map_non_zero} the map
$$ g_i: V_{\lam_i} \lra V_{\lam_0}$$
is non-zero. In particular considering properties of $\gots_{\lam_0-\lam_i} \simeq \gotsl_2$ I have 

$$
g_i(v_i) = g_i(g_{\lam_i-\lam_0}(v_0)) \in V_{\lam_0} \backslash \{0\}
$$  
This in particular means that $k_i$ is non-zero. Therefore:

\begin{cor}\label{corollary_for_the_proposition}
$g_i(v_n)$ is not contained in the tangent space to $\hat{X}$ at the point $v_0$ (i.e. in the
$\sspan \{v_0,\ldots, v_{n-1} \}$). 
\end{cor}
\noprf

\subsection{The semisimple case}\label{semisimple_case}

Just now it is the right time to understand 
 what happens if the group $G$ is not simple, so that $\gotg=\gota \oplus \gotb$ is 
a non trivial splitting. Although this case was dealt with by 
Landsberg and Manivel in \cite[cor.6 and section 2.4]{landsbergmanivel04}, I will provide a different proof. 
Treat this subsection as a digression presenting what kind of a tool am I developing.

\begin{figure}[htb]
\centering
\includegraphics[width=0.55\textwidth]{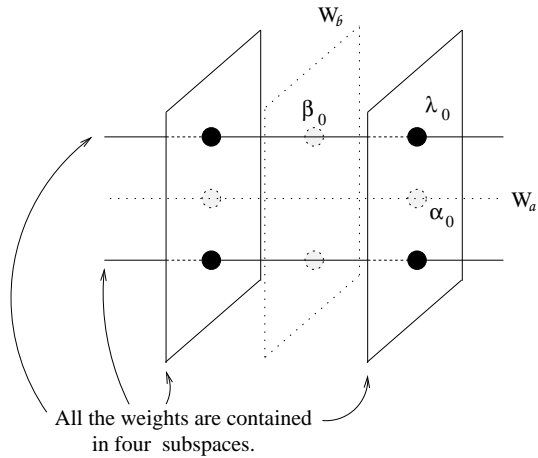}
\caption{The four subspaces with the first coordinate $\pm \alpha_0$  or the second $\pm \beta_0$.}
 \label{pic_semisimple}
\end{figure}

So let $\gotg=\gota \oplus \gotb$ and $V=W_\gota \otimes W_\gotb$ as in the lemma \ref{lemma_for_semisimple}.
Denote by $\alpha_k$'s  all the weights of $W_\gota$ and by $\beta_l$'s all the weights of $W_\gotb$, 
so that the weights of $V$ are exactly $\alpha_k + \beta_l$. On the other hand, each root of $\gotg$ is either a root of
$\gota$ or a root of $\gotb$. This is in particular true for the roots $\lam_0-\lam_i$ (see \ref{positive_roots}). So each $\lam_i$ has either 
the first, or the second coordinate equal to the same coordinate of $\lam_0$. So write $\lam_0 = \alpha_0 + \beta_0$ and
since every weight of $V$ is of the form $\pm \lam_i$ (see \ref{weights_and_root_of_V}(1)),
 it follows that every weight of $V$ is contained in one of the four
following subspaces: either the first coordinate is equal to $\pm \alpha_0$ or the second one is equal to $\pm \beta_0$
 - see the figure \ref{pic_semisimple}.

\begin{rem}\label{remark_faithful_representation}
Both $\alpha_0$ and $\beta_0$ are non-zero.
\end{rem}

\begin{prf}
Indeed, $\lam_0$ is the highest weight, so if one of it's 
coordinates (say, the second one) is zero, then whole the representation has this coordinate equal to zero. But this means,
that whole subalgebra $\gotb$ acts trivially on $V$, contradicting the proposition \ref{kernel_of_action}. 
\end{prf}

\begin{prop}
One of the representations $W_\gota$ or $W_\gotb$ has exactly 2 weights.  
\end{prop}

\begin{prf}
Both representations have at least 2 weights each: $\pm\alpha_0$ and $\pm\beta_0$. If there is another one for each of them,
say $\alpha_1 \ne \pm \alpha_0$ and $\beta_1 \ne \pm\beta_0$, then $\alpha_1 + \beta_1$ is a weight of $V$, but it
 is not contained on any of the four subspaces. 
\end{prf}

Have a quick look at the classification of semisimple Lie algebras (see for example 
\cite[lectures 11-20 and 22]{fultonharris} or \cite[chapter III]{humphreys})
 to see that in fact the Lie algebra (say it is $\gota$) which admits
an irreducible representation with two weights must be isomorphic to $\gotsl_2$: for all the others, the Weyl group doesn't
 preserve any line (neither it can be a direct sum of $\gotsl_2$ and something, via the same argument as in the
 proof of the remark \ref{remark_faithful_representation}).
 
Also $W_\gota$ is two dimensional, so it is standard representation of 
$\gota \simeq \gotsl_2$ and hence $W_{gotb}$ must be of dimension $n$.
Therefore $X$ is a product hypersurface in $\P^1 \times \P^{n-1}$.
It is obvious now, that $X$ must be a product of $\P^1$ and a quadric 
hypersurface in $\P^{n-1}$ (otherwise $X$ is either degenerate or not 
homogeneous).

Hence the theorem:

\begin{theo}
If $\gotg$ is not simple, then $\gotg \simeq \gotsl_2 \oplus \gotso_{n}$ 
and $X$ is a product of line and a quadric hypersurface. 
\end{theo}
\begin{prf}
It follows from the above considerations, that $X$ is isomorphic to $\P^1 \times Q^{n-2}$. 
On the other hand in the section \ref{line_times_quadric} it was computed
explicitly what is the algebra $\gotg$ isomorphic to.
\end{prf}

\subsection{Lengths and angles} \label{section_final_classification}

Recall, that there exists a unique (up to scalar) non-degenerate inner product on the weight lattice invariant under the
 action of the Weyl group (see \cite[\S14.2]{fultonharris})- it is determined by the Killing form. 
Although it is very hard 
to say anything precise about the Killing form on $\gotg$, I will use the form to exclude a lot of the 
configurations. So in this subsection I am going to compare the angles between the roots and their lengths.\\

\begin{rem}\label{remark_no_30_degs}
Note that if $\alpha$ and $\beta$ are roots of $\gotg$, then the angle between them is a multiple of either 
$\frac{\pi}{4}$ or $\frac{\pi}{3}$: it cannot be equal to $\frac{\pi}{6}$ or $\frac{5\pi}{6}$, since the only Lie
algebras, for which these angles occur are products of $\gotg_2$ and something, but this case is excluded by the 
corollary \ref{not_g2} and the analysis of the semisimple case.
In the following analysis I will restrict some rank 2 semisimple subalgebras of $\gotg$ generated by some positive
roots, as in the subsection \ref{subalgebra_s} - I can assume that they are all not of type $G_2$, since the angles
$\frac{\pi}{6}$ and $\frac{5\pi}{6}$ do not appear on in the root system of $\gotg$.
\end{rem}

\begin{figure}[htb]
\centering
\includegraphics[width=1\textwidth]{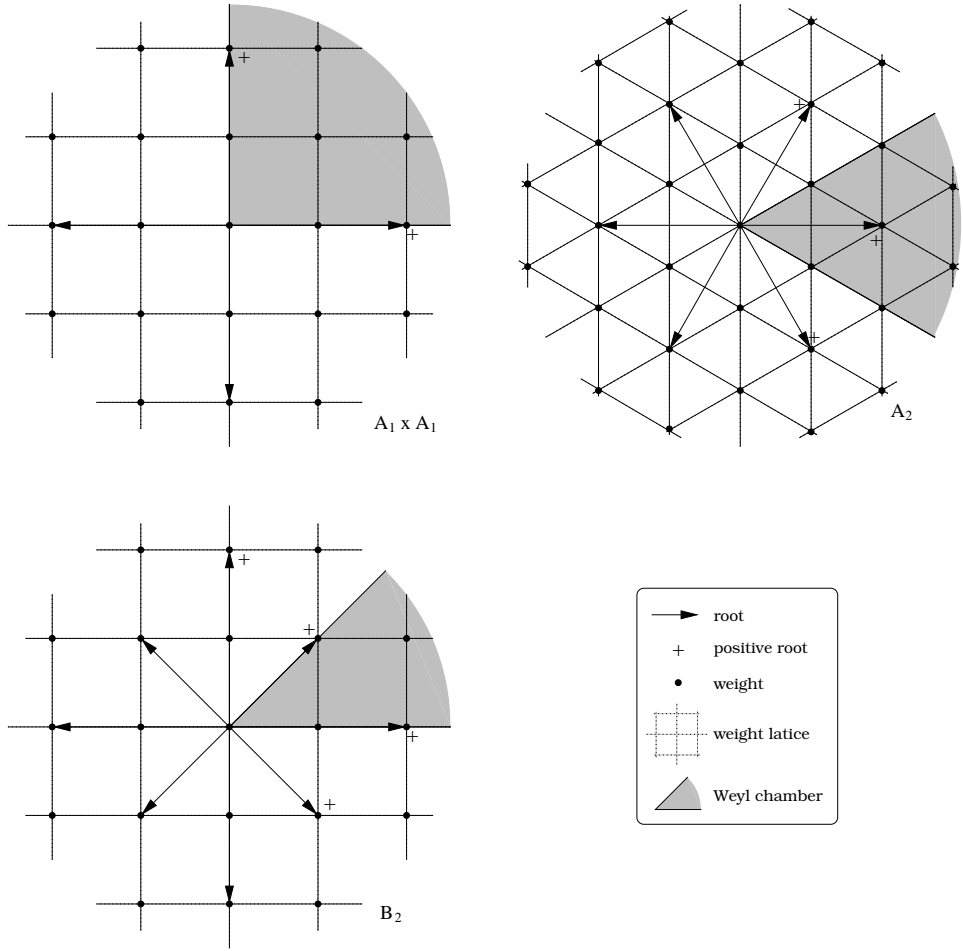}
\caption{The rank two root systems (except $G_2$).}
 \label{pic_root_systems_of_rank_2}
\end{figure}

In the next few proofs I will argue by restricting to some rank 2 subalgebras of $\gotg$. 
So the reader might find it useful to have a quick look at the figure \ref{pic_root_systems_of_rank_2}
to recall the root and weight systems of the rank 2 semisimple Lie algebras.
Note that the system of $G_2$ was omitted since it will not appear in the considerations.

\begin{prop}\label{weights_are_different}
For every $i,j \in \{1,\ldots,n-1\}$ the inequality holds: $\lam_i\ne-\lam_j$.
\end{prop}
\begin{prf}
Assume $\lam_i=-\lam_j$ for some $i,j$. Then in particular $i\ne j$, since $\lam_i \ne 0$.
Moreover both $i,j \ne 0$, since $\lam_0$ is the highest weight and $V$ is irreducible.
Take $\gots<\gotg$ to be the semisimple subalgebra generated by $\alpha_1= \lam_0-\lam_i$ and $\alpha_2=\lam_0-\lam_j$
and the subrepresentation $W \subset V$ with the highest weight $\lam_0$ as in the subsection \ref{subalgebra_s}.
Now $\gots$ is a rank 2 semisimple Lie algebra (not isomorphic to $\gotg_2$ by the remark \ref{remark_no_30_degs}), so it is of type
$A_2$, $B_2$ or $A_1\times A_1$. Note, that by the assumption the weight $\lam_0$ of $\gots$ is a half of a sum of two 
positive roots,
namely: $$\lam_0 = \half (\lam_0-\lam_i \ + \ \lam_0-\lam_j).$$
 Comparing the root systems and weight lattices, this already excludes the case $A_2$. Also in $B_2$ it can only happen if 
$\lam_0-\lam_i$ and $\lam_0-\lam_j$ are
 the shorter perpendicular roots and then their sum (equal to $2\lam_0$) is a root, 
contradicting the corollary \ref{vn_in_X}. 

\begin{figure}[htb]
\centering
\includegraphics[width=0.5\textwidth]{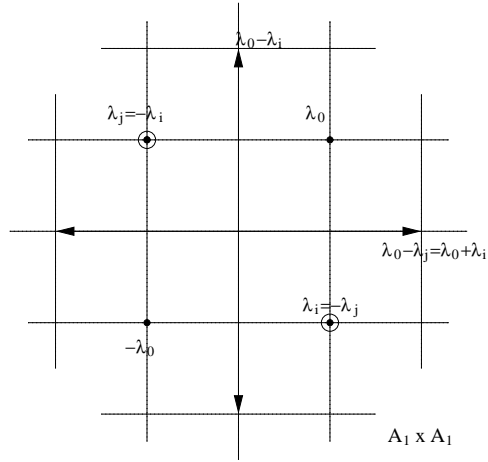}
\caption{The impossible weight diagram for $A_1 \times A_1$ from the proof of the proposition \ref{weights_are_different}.}
 \label{pic_weights_for_A1xA1}
\end{figure}

So the only remaining case is that $\gots$ is of type $A_1 \times A_1$. 
Now since 0 is in the plane passing through the weights $\lam_0, \lam_i, \lam_j$, 
It follows that $-\lam_0$ falls into the slice of weight diagram of $W$ and hence $v_n\in W$.
But then the multiplicity of $\lam_i = -\lam_j$ (and similarly of $\lam_j = -\lam_i$) is equal to 2, since both 
$v_i, g_{\lam_0-\lam_j}(v_n)\in W$ and they are linearly independent by the corollary \ref{corollary_for_the_proposition}.
On the other hand, the multiplicity of $\lam_0$ (and of $-\lam_0$) is 1. Such a weight diagram 
(see the figure \ref{pic_weights_for_A1xA1}) 
cannot occur. 
That excludes the last case and proves the proposition.
\end{prf}

\begin{cor}\label{multiplicity_1}
All the weights of the representation $V$ have multiplicity 1.  
\end{cor}
\begin{prf}
The corollary follows immediately from the corollary \ref{lambdas_are_different}
and the proposition \ref{weights_are_different}.
\end{prf}

\begin{lem}\label{sum_of_2_roots_in_n}
Assume $(\lam_0-\lam_i) + (\lam_0-\lam_j)= (\lam_0-\lam_k)$ for some $i,j,k$. Then the semisimple Lie algebra $\gots$
 generated by any two of them is of type $B_2$ (so $\gots \simeq \gotsp_4$). 
Moreover $(\lam_0-\lam_i)$ and $(\lam_0-\lam_j)$ are perpendicular and 
the representation $W$ as defined in the lemma \ref{subalgebra_s} has exactly 4 weights: $\lam_0, \lam_i, \lam_j,\lam_k$. 
\end{lem} 

\begin{figure}[htb]
\centering
\includegraphics[width=0.5\textwidth]{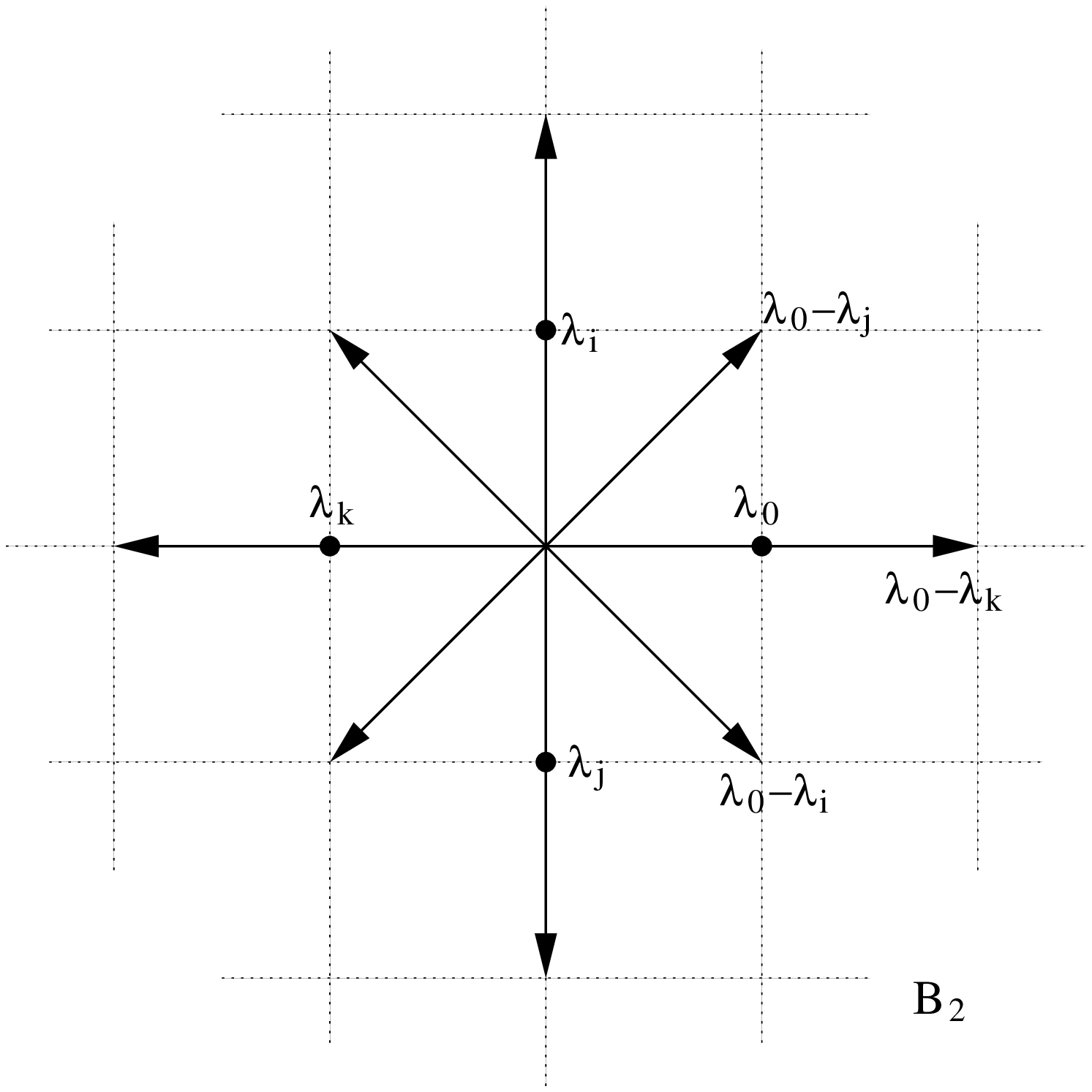}
\caption{The only potentially possible configuration of roots in case
 $(\lam_0-\lam_i) + (\lam_0-\lam_j)= (\lam_0-\lam_k)$ as in the lemma \ref{sum_of_2_roots_in_n}.}
 \label{pic_case_i+j=k}
\end{figure}

\begin{prf}
Note that each pair from $(\lam_0-\lam_i), (\lam_0-\lam_j), (\lam_0-\lam_k)$ is linearly independent, but altogether 
they span
a two dimensional subspace, so $\gots$ does not depend on the choice of the pair. Now $\gots$ is a rank 2 algebra 
(as usually, not isomorphic to $\gotg_2$ by the remark \ref{remark_no_30_degs}), so it is of either of types
 $A_1\times A_1$, $A_2$, $B_2$. The $A_1\times A_1$
 case 
is immediately excluded, since there are only two positive roots there. Now consider the irreducible representation $W$
 generated by $\lam_0$. 

By the corollary \ref{multiplicity_1} all the multiplicities of $W$ are equal to 1. This does not leave a lot of 
possibilities out: in case of $A_2$ (i.e. $\gots=\gotsl_3$), the only such representations are those, whose highest weight 
lays on an edge of the Weyl chamber (see \cite[\S13.2, pp183-184]{fultonharris}). But then only two negative roots move 
the highest weight to another weight and this is not the case for $W$. So indeed the only possibility is that 
$\gots \simeq \gotsp_4$ (so it is of type $B_2$). But then again there are only two possibilities for $W$: it is 
either isomorphic to the natural representation of $\gotsp_4$ or to the subrepresentation of the exterior square of 
the natural one (see  \cite[\S16.2]{fultonharris}). The latter case is in contradiction to  
$(\lam_0-\lam_i) + (\lam_0-\lam_j)= (\lam_0-\lam_k)$. 
So $W$ is the natural representation of $\gotsp_4$ as pictured on the figure \ref{pic_case_i+j=k}.
And this implies the lemma.
\end{prf}

\begin{rem}
Actually, it never happens that $(\lam_0-\lam_i) + (\lam_0-\lam_j)= (\lam_0-\lam_k)$.
\end{rem}

Let $\gotn<\gotg$ be the nilpotent subalgebra generated by all $\gotg_{\lam_i-\lam_0}$'s,
 so that $\gotg= \gotp \oplus \gotn$
 and $R= R_\gotp \cup R_\gotn$, where $R_\gotp$ and $R_\gotn$ are 
the sets of roots of $\gotp$ and $\gotn$ correspondingly.
Now the above lemma says a lot about the configuration of the roots in $R_\gotn$:
\begin{cor}\label{angle_not_obtuse}
Let $\lam_i-\lam_0$ and $\lam_j-\lam_0$ be any two distinct roots from $R_\gotn$. 
Then the angle between the two roots is either 
acute or right. 
\end{cor}
\begin{prf}
If the angle is obtuse, then the sum of the two roots is a root (see \cite[\S21.1, property (6), p324]{fultonharris}) and 
it is in $R_\gotn$, so it is of the form $\lam_k - \lam_0$, contradicting the lemma.
\end{prf}

\begin{theo}\label{one_simple_root}
If $\gotg$ is simple then only one simple root belongs to $R_\gotn$.
\end{theo}
\begin{prf}
This is a consequence of the proposition \ref{one_simple_root_generally} and the corollary \ref{angle_not_obtuse}. 
\end{prf}

Note that the theorem already implies the classification: if only one simple root belongs to $R_\gotn$, then the Picard
number of $X$ is one ([citation!!!])
 and the homogeneous legendrian varieties with Picard number one are only the subadjoint 
varieties (see \cite[\S2.4, thm 11]{landsbergmanivel04}). Yet I will proceed with the proof to provide another
approach.

\subsection{Classification}

Let me summarise, what properties of the irreducible representation $V$ are proved so far:
\begin{itemize}
\item[(i)] 
The highest weight of $V$ lies on an edge of the Weyl chamber 
(this follows from the theorem \ref{one_simple_root}).
\item[(ii)] 
The dimension of $V$ is twice the dimension of the orbit of the highest weight vector $v_0$ (i.e. the number of roots
 moving the highest weight to some other weight plus one).
\item[(iii)]
$V$ is isomorphic to $V^*$.
\item[(iv)] All the weights have multiplicity 1 (see the corollary \ref{multiplicity_1}).
\end{itemize}

The first two properties are the most restrictive. So the strategy now is:
\begin{itemize}
\item[1)]
Fix a simple group.
\item[2)]
Choose an edge of it's Weyl Chamber.
\item[3)]
Ask whether the representation with the highest weight on the edge satisfies (ii)-(iv):
\begin{itemize}
\item[a)]
If the dimension of the representation is greater than twice the dimension of the orbit then no representation from
this edge can satisfy (ii) (all the further weights correspond to even bigger representations). 
\item[b)]
If the first dimension is less that the second, then try with the next weight along the same edge.
 \end{itemize}
\item[4)]
Finally if the dimensions satisfy (ii), check whether (iii) and (iv) are true.
\end{itemize}
After running through this algorithm one is left only with very few possibilities. Let me  see how does it work for 
classical Lie algebras:

\begin{theo}
If $\gotg$ is of type $A_m$ for some $m$, then either $m=1$ and $X$ is the twisted cubic or $m=5$ and 
$X$ is the Grassmannian $Gr(3,6)$. 
\end{theo} 
\begin{prf}
For this case the most restrictive are the properties (i) and (iii) -
together they imply that $V$ is in fact contained as the biggest subrepresentation in some symmetric power
of \textbf{the middle exterior power} of the natural representation and hence $m$ is odd and 
$X \simeq Gr(\frac{m+1}{2}, m+1)$. For the dimension reason (see 3a) of the algorithm)
 $m$ cannot be greater than $5$. 

For $m=1$ one gets the twisted cubic. 

For $m=3$ one can see that if $\lam_0$ is the first weight along the edge, than the dimension of $V$ would be
$6$ while the dimension of $\hat{X}$ would be $5 \ne \half \cdot 6$. On the other hand if $\lam_0$ is a further weight
than the first one, then the dimension of $V$ would be greater than $20 > 2 \cdot 5$. 

For $m=5$ one gets the Grassmannian $Gr(3,6)$ which is a legendrian variety under the Pl\"ucker embedding 
- see \ref{example_grassmannian}.  
\end{prf}

Similar proof is going to work for the $C_m$ case (including $C_2=B_2$):
\begin{theo}
If $\gotg$ is of type $C_m$ for some $m$, then $m=3$ and $X$  is the Lagrangian Grassmannian $G_L(3,6)$. 
\end{theo} 
\begin{prf}
The theorem \ref{one_simple_root} implies that $X\simeq Gr_L(k,2m)$ for some $k \le m$. 
Applying 3a) one gets 
 that either $k=1$ (so $X\simeq \P^{2i-1}$) or $k=2, m=2$ or  $k=2, m=3$ or $k=3, m=3$.

To exclude the case $k=1$, let $\alpha$ be the highest weight of the natural representation of $\gotsp_{2m}$. Then $\lam_0$
must be a multiple of $\alpha$. It cannot be just equal to $\alpha$ (because $X \ne \P(V)$). But all the further 
representations contain a multiple weight, contradicting (iv). 

Both $k=2$ cases are excluded by computing the dimensions of the proper representations - 3b) and 3a).

Finally $k=3, m=3$ is the required case: $X\simeq Gr_L(3,6)$ and it is a legendrian variety
 - see  \ref{example_lagrangian_grassmannian}.

\end{prf}

Analogous argument works for also for $B_{m_1}$ and $D_{m_2}$:
\begin{theo}
If $\gotg$ is isomorphic to $\gotso_m$ for some $m \ge 7$, then 
$m=12$ and $X$ is the spinor variety $\mathbb{S}_6$.
\end{theo}
\begin{prf}
Argue in the same way as previously: first compare the dimensions of the minimal representations and the dimensions of
the corresponding homogeneous spaces (which  in these case are the Grassmannians $Gr_o(k,m)$ of 
k-hyperplanes isotropic with respect to a nondegenerate quadratic form on $\C^m$). Conclude that the dimension of the 
representation is to big, unless $k=1$ (so that $X \simeq Gr_o(1,2m+1) \simeq Q^{2m-1}$) or $k=m \in \{7,8,9,10,11,12\}$
 (so that $X$ is the maximal isotropic Grassmannian, so called spinor variety denoted 
$\mathbb{S}_{[\frac{m}{2}]}$).

The $k=1$ case and $k=m \in \{7,8,9,10\}$ are excluded exactly in the same way as for $C_m$.

As for the $k=m=11$ or $k=m=12$, the spinor variety is the same for both cases and it is a legendrian variety.
As it has been stated in the subsection  \ref{example_spinor}, the Lie algebra $\gotg$ is of type $D_6$ ($\gotso_{12}$)
and  hence $k=m=12$.
\end{prf}

Now the only remaining cases are that $\gotg$ is one of the exceptional algebras: $F_4, E_6,E_7,E_8$ (the $G_2$ case was 
excluded in corollary \ref{not_g2}). Since in these cases the computations of the dimensions and multiplicities are 
rather complicated to deal with, I refer to computer calculation,
(which implements the algorithm above - there are only finite number of things to check) to state that:

\begin{theo}
If $\gotg$ is a exceptional Lie algebra, then it is of type $E_7$ and $X$ is a 27-dimensional homogeneous variety
embedded in $\P^{55}$ - see \ref{example_E7}.
\end{theo}
\noprf

Summarizing the above theorems one gets the list of subadjoint varieties:

\begin{theo}\label{final_classification}
If $X$ is a smooth irreducible legendrian variety generated by quadrics, then $X$ is one of the subadjoint varieties
(see the section \ref{rozdzial_przyklad} and in particular the table \ref{table_contacts}
 on page \pageref{table_contacts}):
\begin{itemize}
\item the twisted cubic: $\P^1 \subset \P^3$;
\item the product of a line and a quadric: $\P^1 \times Q^{n-2} \subset \P^{2n-1}$; 
\item the Grassmannian of Lagrangian subspaces in $\C^6$: $Gr_{Lag}(3,6) \subset \P^{13}$;
\item the Grassmannian $Gr(3,6) \subset \P^{19}$;
\item the spinor variety $\mathbb{S}_6 \subset \P^{31}$;
\item the $E_7$ variety  $X^{27} \subset \P^{55}$.
\end{itemize}
\end{theo}

\end{document}